\documentclass[11pt]{amsart}
\usepackage{amsmath,amscd,amssymb,amsfonts,latexsym}
\usepackage[all]{xy}
\usepackage{ulem}

\usepackage[latin1]{inputenc}

\usepackage[T1]{fontenc}

\usepackage{eucal} 
\usepackage{cmmib57} 
\usepackage{amsbsy}
\usepackage{amsthm,alltt,dsfont}

\usepackage{geometry}
\allowdisplaybreaks

\newcommand{\bC}{{\mathbb C}}

\newcommand{\bN}{{\mathbb N}}

\newcommand{\bQ}{{\mathbb Q}}
\newcommand{\bZ}{{\mathbb Z}}

\newcommand{\cD}{{\mathcal D}}

\newcommand{\cF}{{\mathcal F}}

\newcommand{\cM}{{\mathcal M}}

\newcommand{\cL}{{\mathcal L}}
\newcommand{\cV}{{\mathcal V}}

\newcommand{\Zs}{{Z^{(n)}}}

\newcommand{\co}{\it Coh}
\newcommand{\cons}{\it Constr}
\newcommand{\fo}{\it For}

\newcommand{\mhs}{{\it MHS}^p}

\newcommand{\lra}{\longrightarrow}

\newcommand{\MHM}{{\it MHM}}

\theoremstyle{plain}
\newtheorem{thm}{Theorem}[section]
\newtheorem{cor}[thm]{Corollary}

\newtheorem{lem}[thm]{Lemma}
\newtheorem{prop}[thm]{Proposition}

\newtheorem*{ack}{Acknowledgements}

\theoremstyle{definition}
\newtheorem{df}[thm]{Definition}
\newtheorem{rem}[thm]{Remark}
\newtheorem{example}[thm]{Example}

\def\be{\begin{equation}}
\def\ee{\end{equation}}

\def\bt{\begin{thm}}
\def\et{\end{thm}}

\def\bc{\begin{cor}}
\def\ec{\end{cor}}

\def\br{\begin{rem}}
\def\er{\end{rem}}

\def\bp{\begin{prop}}
\def\ep{\end{prop}}

\def\bl{\begin{lem}}
\def\el{\end{lem}}

\def\bn{\begin{enumerate}}
\def\en{\end{enumerate}}

\def\bex{\begin{example}}
\def\eex{\end{example}}

\def\bd{\begin{df}}
\def\ed{\end{df}}

\title{Equivariant characteristic classes of external and symmetric products of varieties}
\author{Lauren\c{t}iu Maxim}
\address{L. Maxim: Department of Mathematics, University of Wisconsin-Madison, 480 Lincoln Drive, Madison, WI 53706, USA}
\email {maxim@math.wisc.edu}
\author[J. Sch\"urmann ]{J\"org Sch\"urmann}
\address{J.  Sch\"urmann : Mathematische Institut,
          Universit\"at M\"unster,
          Einsteinstr. 62, 48149 M\"unster,
          Germany.}
\email {jschuerm@uni-muenster.de}

\date{\today}



\begin{document}

\begin{abstract} We obtain refined generating series formulae for equivariant characteristic classes of external and symmetric products of singular complex quasi-projective varieties. More concretely, we study equivariant versions of Todd, Chern and Hirzebruch classes for singular spaces, with values in delocalized Borel-Moore homology of external and symmetric products. As a byproduct, we recover our previous characteristic class formulae for symmetric products and obtain new equivariant generalizations of these results, in particular also in the context of twisting by representations of the symmetric group. \end{abstract}

\maketitle

\tableofcontents

\section{Introduction}
In this paper, we obtain refined generating series formulae for equivariant characteristic classes of external and symmetric products of singular complex quasi-projective varieties, generalizing our previous results for symmetric products from \cite{CMSSY}.

\subsection{Equivariant characteristic classes} All spaces in this paper are assumed to be complex quasi-projective, though many constructions also apply to other categories of spaces with a finite group action (e.g., compact complex analytic  manifolds or varieties over any base field of characteristic zero). 
For such a variety $X$, consider an algebraic action $G\times X \to X$ by a finite group $G$, with quotient map $\pi:X \to X':=X/G$. For any  $g \in G$, we let $X^g$ denote the corresponding  fixed point set. 

We let $cat^G(X)$ be a category of $G$-equivariant objects on $X$ in the underlying category $cat(X)$ (e.g., see \cite{CMSS,MS}), which in this paper refers to one of the following examples: coherent sheaves $\co(X)$, algebraically constructible sheaves of complex vector spaces $\cons(X)$, and (algebraic) mixed Hodge modules $\MHM(X)$ on $X$. 
We denote by $K_0(cat^G(X))$ the corresponding Grothendieck groups of these $\bQ$-linear abelian categories. We will also work with the relative Grothendieck group $K_0^G(var/X)$ of $G$-equivariant quasi-projective varieties over $X$, defined by using the scissor relation as in \cite{CMSS}.  
Let  $H_*(X)$ denote the even degree Borel-Moore homology $H^{BM}_{ev}(X) \otimes R$ with coefficients in a commutative $\bC$-algebra $R$ (resp. $\bQ$-algebra if $G$ is a symmetric group). Note that $H_*(-)$ is functorial for all proper maps, with a compatible cross-product $\boxtimes$. 

Let 
$$cl_*(-;g):K_0(cat^G(X)) \to H_{*}(X^g)$$
be one of the following {\it equivariant characteristic class transformation of Lefschetz type}, see Section \ref{s3b}:
\begin{enumerate}
\item the Lefschetz-Riemann-Roch transformation of Baum-Fulton-Quart \cite{BFQ} and Moonen \cite{Mo}:
$$
td_*(-;g):K_0(\co^G(X)) \lra H_{*}(X^g),
$$
with $R=\bC$ (resp. $R=\bQ$  if $G$ is a symmetric group).
\item the localized Chern class transformation \cite{Sch1}:
$$c_*(-;g):K_0(\cons^G(X)) \lra H_{*}(X^g),$$
with $R=\bC$ (resp. $R=\bQ$  if $G$ is a symmetric group).
\item the motivic version of the (un-normalized) Atiyah-Singer  class transformation \cite{CMSS}:
$$T_{y*}(-;g):K^G_0 (var/X) \lra H_{*}(X^g),$$
with $R=\bC[y]$ (resp. $R=\bQ[y]$  if $G$ is a symmetric group).
\item the mixed Hodge module version of the (un-normalized) Atiyah-Singer  class transformation \cite{CMSS}:
$$T_{y*}(-;g): K_0 (\MHM^G(X)) \lra H_{*}(X^g),$$
with $R=\bC[y^{\pm 1}]$ (resp. $R=\bQ[y^{\pm 1}]$  if $G$ is a symmetric group).
\end{enumerate}
These class transformations are covariant functorial for $G$-equivariant proper maps and cross-products $\boxtimes$. Over a point space, they reduce to a certain $g$-trace (as explained in Section \ref{s3b}). For a subgroup $K$ of $G$, with $g \in K$, these transformations $cl_*(-;g)$ of Lefschetz type commute with the obvious restriction functor ${\rm Res}^G_K$.  Moreover, $cl_*(-;g)$ depends only on the action of the cyclic subgroup generated by $g$. In particular, if $g=id_G$ is the identity of $G$, we can take $K$ the identity subgroup $\{id_G\}$ with  ${\rm Res}^G_K$ the forgetful functor $${\it For}: K_0(cat^G(X)) \to K_0(cat(X)),$$ so that $cl_*(-;id_G)=cl_*(-)$ fits with a corresponding non-equivariant characteristic class, which in the above examples are: 
\begin{enumerate}
\item the Todd class transformation $td_*$ of Baum-Fulton-MacPherson \cite{BFM} appearing in  the Riemann-Roch theorem for singular varieties.
\item the MacPherson-Chern class transformation $c_*$ \cite{MP}.
\item the motivic version of the (un-normalized) Hirzebruch class transformation  $T_{y*}$ of \cite{BSY}.
\item the mixed Hodge module version of the (un-normalized) Hirzebruch class transformation  $T_{y*}$ of \cite{BSY}, see also \cite{Sch2}.
\end{enumerate}

\medskip

The disjoint union $IX:=\bigsqcup_{g \in G} X^g$ (which is also called the inertia space of the $G$-space $X$) admits an induced $G$-action by $h: X^g \to X^{hgh^{-1}}$, such that the canonical map $$i: IX=\bigsqcup_{g \in G} X^g \to X$$ defined by the inclusions of fixed point sets becomes $G$-equivariant. Therefore, $G$ acts in a natural way on  $\bigoplus_{g \in G} \ H_*(X^g)=H_*(IX)$ by conjugation.
\bd The {\it (delocalized) $G$-equivariant homology of $X$} is the $G$-invariant subgroup \be H_*^{G}(X):=\big( H_*(IX) \big)^G=\left( \bigoplus_{g \in G} \ H_*(X^g)\right)^G.\ee
This theory is functorial for proper $G$-maps and induced cross-products $\boxtimes$.
 \ed
 This notion is different (except for free actions) from the equivariant Borel-Moore homology $H^G_{BM,2*}(X)\otimes R$ defined by the Borel construction. In fact, since $G$ is finite and $R$ is a $\bQ$-algebra, one has \be\label{idnew} H^G_{BM,2*}(X) \otimes R\simeq \left(H^{BM}_{2*}(X)\otimes R  \right)^G \simeq H^{BM}_{2*}(X/G)\otimes R,\ee which is just a direct summand  of $H^{G}_{*}(X)$ corresponding to the identity element of $G$, denoted by $H^G_{id,*}(X)$.
For example, if $G$ acts {\it trivially} on $X$ (e.g., $X$ is a point), then 
 $$ H_*^{G}(X) \simeq H_*(X) \otimes  C(G),$$
 where $C(G)$ denotes the free abelian group of $\bZ$-valued  class functions on $G$ (i.e., functions which are constant on the conjugacy classes of $G$). Note also that 
 \be\label{bc} H_*^{G}(X)=\big( H_*(IX) \big)^G \simeq H_*(IX/G).\ee
 
\bd For any of the Lefschetz-type characteristic class transformations ${cl}_*(-;g)$ considered above, we define a corresponding  $G$-equivariant  class transformation (with $T^G_{y*}$ the $G$-equivariant Hirzebruch class transformation)
$$cl^G_{*}: K_0(cat^G(X)) \to  H^{G}_{*}(X)$$
by: $$cl^G_{*}(-):=\bigoplus_{g \in G} cl_*(-;g) \in \left( \bigoplus_{g \in G} \ H_*(X^g) \right)^G.$$
\ed
\noindent The $G$-invariance of the class $cl_*^G(-)$ is a consequence of {\it conjugacy invariance} of  the Lefschetz-type characteristic class ${cl}_*(-;g)$, see \cite{CMSS}[Sect.5.3]. Note that the summand $cl_*(-;id)\in (H_*(X))^G$ corresponding to the identity element of $G$ is just the non-equivariant characteristic class, which for equivariant coefficients is invariant under the $G$-action by functoriality. Under the identification (\ref{idnew}), this class also agrees (for our finite group $G$) with the corresponding (naive) equivariant characteristic class defined in terms of the Borel construction, e.g., for $cl_*=td_*$, this is the equivariant Riemann-Roch-Transformation of Edidin-Graham \cite{EG}; and for $cl_*=c_*$, this is the equivariant Chern class transformation of Ohmoto \cite{Oh1,Oh2}.

The above transformation $cl^G_*(-)$ has the same properties as the Lefschetz-type transformations $cl_*(-;g)$, e.g., functoriality for proper push-downs, restrictions to subgroups, and multiplicativity for exterior products. 

\br The $G$-equivariant characteristic classes defined here for $cl_*=td_*, T_{y*}$ agree, up to the normalization factor $\frac{1}{|G|}$, with the corresponding notions introduced in \cite{CMSS}.
\er


\subsection{Generating series formulae}\label{objects}

Let now $Z$ be a quasi-projective variety, and denote by $\Zs:=Z^n/{\Sigma_n}$ its $n$-th symmetric product (i.e., the quotient of $Z^n$ by the natural permutation action of the symmetric group $\Sigma_n$ on $n$ elements), with $\pi_n:Z^n \to \Zs$ the natural projection map. The standard approach for computing invariants of the symmetric products $\Zs$ is to collect the respective invariants of all symmetric products in a generating series, and then compute the latter solely in terms of invariants of $Z$, e.g., see \cite{CMSSY} and the references therein. In this paper, we obtain generalizations of results of \cite{CMSSY}, formulated in terms of equivariant characteristic classes of external products and resp., symmetric products of varieties.

To a given object $\cF \in cat(Z)$ in a category as above, i.e., coherent or constructible sheaves, or mixed Hodge modules on $Z$ (resp., morphisms $f:Y\to Z$ in the motivic context), we attach new objects  as follows (see \cite{CMSSY,MSS,MS} for details):
\bn
\item[(a)] the $\Sigma_n$-equivariant object $\cF^{\boxtimes n} \in cat^{\Sigma_n}(Z^n)$ on the cartesian product $Z^n$ (e.g., $f^n:Y^n \to Z^n$ in the motivic context).
\item[(b)] the $\Sigma_n$-equivariant object $\pi_{n*}\cF^{\boxtimes n}\in cat^{\Sigma_n}(\Zs)$ on the symmetric product $\Zs$ (e.g., the $\Sigma_n$-equivariant map $Y^n \to \Zs$ in the motivic context).
\item[({c})] the following non-equivariant objects in $cat(\Zs)$:
\begin{enumerate}
\item[(1)] the $n$-th {\it symmetric power} object $\cF^{(n)}:=\left( \pi_{n*}\cF^{\boxtimes n} \right)^{\Sigma_n}$ on $\Zs$, defined by using the projector $(-)^{\Sigma_n}$ onto the $\Sigma_n$-invariant part (respectively, the map $f^{(n)}:Y^{(n)} \to Z^{(n)}$ induced by dividing out the $\Sigma_n$-action in the motivic context).
\item[(2)] the $n$-th {\it alternating power} object $\cF^{\{n\}}:=\left( \pi_{n*}\cF^{\boxtimes n} \right)^{sign-\Sigma_n}$ on $\Zs$, defined by using the alternating projector $(-)^{sign-\Sigma_n}$ onto the sign-invariant part.
 (This construction does not apply in the motivic context.)
\item[(3)]  $\fo(\pi_{n*}\cF^{\boxtimes n})$, obtained by forgetting the $\Sigma_n$-action on $\pi_{n*}\cF^{\boxtimes n} \in cat^{\Sigma_n}(\Zs)$  (e.g., the induced map $Y^n \to \Zs$ in the motivic context).
 \end{enumerate}
\en
These constructions and all of the following results also apply to suitable bounded complexes (e.g., the constant Hodge module complex $\bQ^H_Z$), see Remark \ref{5.9} for details. 

\medskip

The main goal of this paper is to compute generating series formulae for the (equivariant) characteristic classes of these new coefficients only in terms of the original characteristic class $cl_*(\cF)$. These generating series take values in a corresponding commutative graded $\bQ$-algebra $\mathds{H}^{\Sigma}_*(Z)$, $\mathds{PH}^{\Sigma}_*(Z)$ and resp. $\mathds{PH}_*(Z)$, and are formulated with the help of certain operators which transport homology classes from $Z$ into these corresponding commutative graded $\bQ$-algebras. In each of three situations (a)-(c) above, these algebras of Pontrjagin type and operators are described explicitely as follows:
\bn
\item[(a)] $$\mathds{H}^{\Sigma}_*(Z):=\bigoplus_{n \geq 0} H_*^{\Sigma_n}(Z^n) \cdot t^n, $$ with {\it creation operator} $\mathfrak{a}_r$ defined by: if $\sigma_r=({r})$ is an $r$-cycle in $\Sigma_r$, then $\mathfrak{a}_r$ is the composition:
$$\mathfrak{a}_r:H_*(Z) \overset{\cdot r}{\to} H_*(Z) \cong H_*((Z^r)^{\sigma_r})^{Z_{\Sigma_r}(\sigma_r)} \hookrightarrow H_*^{\Sigma_r}(Z^r),$$
where $\langle \sigma_r \rangle =Z_{\Sigma_r}(\sigma_r)$ acts trivially on $H_*((Z^r)^{\sigma_r})$.

\noindent Note that the direct summand $$\mathds{H}^{\Sigma}_{id,*}(Z):=\bigoplus_{n \geq 0} H_{id,*}^{\Sigma_n}(Z^n) \cdot t^n \subset \mathds{H}^{\Sigma}_*(Z)$$ corresponding to the identity component is a subring, and the projection of $\mathds{H}^{\Sigma}_*(Z)$ onto the subring $\mathds{H}^{\Sigma}_{id,*}(Z)$ kills all the creation operators except $\mathfrak{a}_1=id_{H_*(Z)}.$

\item[({b})]  \begin{align*} \mathds{PH}_*^{\Sigma}(Z)&:=\bigoplus_{n \geq 0} H_*^{\Sigma_n}(\Zs) \cdot t^n \\
&\simeq \bigoplus_{n \geq 0} \left( H_*(\Zs) \otimes C(\Sigma_n) \right) \cdot t^n \hookrightarrow   \mathds{PH}_*(Z) \otimes \bQ[p_i \ \vert \ i \geq 1],\end{align*}
with corresponding operator 
$p_r \cdot d_{r*}:H_*(Z) \to H_*(Z^{(r)}) \otimes \bQ[p_i, i \geq 1]$, where $d_r:=\pi_r \circ \Delta_r:Z \to Z^{({r})}$ is the composition of the natural projection $\pi_r:Z^r \to Z^{({r})}$ with the diagonal embedding $\Delta_r:Z \to Z^r$. The algebra inclusion above is induced from the Frobenius character 
$$ch_F: C(\Sigma)\otimes \bQ:=\oplus_n C(\Sigma_n) \otimes \bQ \overset{\simeq}{\lra} \bQ[p_i, i \geq 1]=:\Lambda \otimes \bQ$$
to the graded ring of $\bQ$-valued symmetric functions in infinitely many variables $x_m$ ($m \in \bN$), with $p_i:=\sum_m x_m^i$ the $i$-th power sum function, see \cite{Mc}[Ch.I,Sect.7], and $\mathds{PH}_*(Z)$ defined as below. 

\item[({c})] $$\mathds{PH}_*(Z):=\bigoplus_{n \geq 0} H_*(\Zs) \cdot t^n,$$ with corresponding operator $d_{r*}:H_*(Z) \to H_*(Z^{(r)})$.
\en
Moreover, the creation operator $\mathfrak{a}_r$ satisfies the identity:
$$ \pi_{r*} \circ \mathfrak{a}_r=p_r \cdot d_{r*},$$ justifying the multiplication by $r$ in its definition.

\bigskip

The main characteristic class formula of this paper is contained in:  
\bt\label{thint} The following generating series formula holds in the commutative graded $\bQ$-algebra  $\mathds{H}^{\Sigma}_*(Z):=\bigoplus_{n \geq 0} H_*^{\Sigma_n}(Z^n)\cdot t^n$, if $cl_*$ is either $td_*$, $c_*$ or $T_{-y*}$:
\be\label{m2i} \sum_{n \geq 0} cl_*^{\Sigma_n} (\cF^{\boxtimes n}) \cdot t^n=\exp \left( \sum_{r\geq 1} \mathfrak{a}_r (\Psi_r(cl_*(\cF))\cdot \frac{t^r}{r} \right),\ee
where $\Psi_r$ denotes the homological Adams operation defined by
$$\Psi_r=
\begin{cases}
id & {\rm if} \ cl_*=c_*\\
\cdot \frac{1}{r^i} \ {\rm on} \ H^{BM}_{2i}(Z) \otimes \bQ & {\rm if} \ cl_*=td_* \\
\cdot \frac{1}{r^i} \ {\rm on} \ H^{BM}_{2i}(Z) \otimes \bQ, {\rm and} \ y \mapsto y^r & {\rm if} \ cl_*=T_{-y*}.
\end{cases}
$$
In particular, by projecting onto the identity component, we get
\be\label{m2j} \sum_{n \geq 0} cl_* (\cF^{\boxtimes n};id) \cdot t^n=\exp \left( t \cdot cl_*(\cF) \right) \in \mathds{H}^{\Sigma}_{id,*}(Z).
\ee
\et
For the rest of this Introduction, $cl_*$ denotes any of the classes $td_*$, $c_*$ or $T_{-y*}$.
The proof of Theorem \ref{thint} is purely formal, based on the multiplicativity and conjugacy invariance of the Lefschetz-type claracteristic classes $cl_*(-;g)$, together with the following key {\it localization} formula from \cite{CMSSY}[Lemma 3.3, Lemma 3.6, Lemma 3.10]:
\be\label{kl}
cl_*(\cF^{\boxtimes r}; \sigma_r)=\Psi_r cl_*(\cF), 
\ee
under the identication $(Z^r)^{\sigma_r} \simeq Z$. For this localization formula in the context of Hirzebruch classes, it is important to work with the parameter $-y$ and the un-normalized versions of Hirzebruch classes and their respective equivariant analogues, see \cite{CMSSY}.
In fact, formula (\ref{m2i}) is a special case of an abstract generating series formula (\ref{m2}), which
holds for any functor $H$ (covariant for isomorphisms)
with a compatible commutative, associative cross-product $\boxtimes$, with a unit $1_{pt}\in H(pt)$. The above-mentioned abstract formula (\ref{m2}) codifies the combinatorics of the action of the symmetric groups $\Sigma_n$,
and it should be regarded as a far-reaching generalization of the well-known identity of symmetric functions (e.g., see the proof of \cite{Mc}[(2.14)]):
$$ \sum_{n \geq 0}  h_n t^n=\exp \left( \sum_{r\geq 1}  p_r \cdot \frac{t^r}{r} \right),
$$
with $h_n$ the $n$-th complete symmetric function.
Other applications of the abstract generating series formula (\ref{m2}) in the framework of orbifold cohomology and, resp., localized $K$-theory are explained in Section \ref{exabs}. In this way, we reprove and generalize some results from \cite{QW} and, resp., \cite{W}. Moreover, in Section \ref{genser}, we give another application of (\ref{m2}) to canonical constructible functions and orbifold-type Chern classes of symmetric products, reproving some results of Ohmoto \cite{Oh2}.

\br
In the motivic context, the exponentiation map
\be\label{moti} K_0(var/Z) \lra \bigoplus_{n\geq 0} K_0^{\Sigma_n}(var/Z^n) \cdot t^n, \ [f:X \to Z] \mapsto \sum_{n \geq 0} [f^n:X^n \to Z^n] \cdot t^n
\ee
is well-defined as in \cite{Be}, and should be regarded as an equivariant analogue of the (relative) Kapranov zeta function used in \cite{CMSSY,MS}. In fact, the latter can be recovered from (\ref{moti}) by pushing down to the symmetric products (resp. to a point), and taking the quotients by the $\Sigma_n$-action.
\er

\medskip

By pushing formula (\ref{m2i}) down to the symmetric products, we obtain by functorialy the following result:
\bc\label{cormaini}
The following generating series formula holds in the commutative graded $\bQ$-algebra $\mathds{PH}^{\Sigma}_*(Z):= \bigoplus_{n \geq 0} H_*^{\Sigma_n}(\Zs) \cdot t^n \hookrightarrow \mathds{PH}_*(Z) \otimes \bQ[p_i , i \geq 1]$:
\be\label{ci}  \sum_{n \geq 0} cl_*^{\Sigma_n} (\pi_{n*}\cF^{\boxtimes n}) \cdot t^n=\exp \left( \sum_{r\geq 1} p_r\cdot d_{r*} (\psi_r(cl_*(\cF))\cdot \frac{t^r}{r} \right).
\ee
\ec
This should be regarded as a characteristic class version of \cite{Ge}[Prop.5.4].
In particular, if $Z$ is projective, then by taking degrees, we get in Section \ref{gensersp} generating series formulae for the characters of virtual $\Sigma_n$-representations of $H^*(Z^n;\cF^{\boxtimes n})$, that is, 
\be\label{deg1i}  \sum_{n \geq 0} tr_{\Sigma_n} (Z^n;\cF^{\boxtimes n}) \cdot t^n=\exp \left( \sum_{r\geq 1} p_r\cdot \chi(H^*(Z,\cF)) \cdot \frac{t^r}{r} \right) \in \bQ[p_i , i \geq 1][[t]],
\ee
for $\cF$ a coherent or constructible sheaf and with $\chi$ denoting the corresponding Euler characteristic, respectively,
\be\label{deg2i}  \sum_{n \geq 0} tr_{\Sigma_n} (Z^n;\cM^{\boxtimes n}) \cdot t^n=\exp \left( \sum_{r\geq 1} p_r\cdot \chi_{-y^r}(H^*(Z,\cM)) \cdot \frac{t^r}{r} \right) \in \bQ[y^{\pm 1}, p_i , i \geq 1][[t]], 
\ee
for $\cM$ a mixed Hodge module on $Z$, and with $\chi_y(H^*(Z,\cM))$ the corresponding $\chi_y$-polynomial.

\medskip

By specializing all the $p_i$'s to the value $1$ (which corresponds to the use of the projectors $(-)^{\Sigma_n}$), formula (\ref{ci}) reduces to 
 the main result of \cite{CMSSY}, namely:
\bc\label{c1.6} The following generating series formula holds in the Pontrjagin ring $\mathds{PH}_*(Z):=\bigoplus_{n \geq 0} H_*(\Zs) \cdot t^n$:
\be \sum_{n \geq 0} cl_* (\cF^{(n)}) \cdot t^n=\exp \left( \sum_{r\geq 1}  d_{r*} (\psi_r(cl_*(\cF))\cdot \frac{t^r}{r} \right).
\ee
\ec
In particular, if $Z$ is projective, we recover the degree formulae from \cite{CMSSY}, which can now be also derived from (\ref{deg1i}) and (\ref{deg2i}) by specializing all $p_i$'s to $1$.

\medskip

Corollary \ref{cormaini} also has other important applications.
For example, by specializing the $p_i$'s to the value $sign(\sigma_i)=(-1)^{i-1}$ (which corresponds to the use of the alternating projectors $(-)^{sign-\Sigma_n}$), formula (\ref{ci}) reduces to:
\bc\label{c1.7} The following generating series formula holds in the Pontrjagin ring $\mathds{PH}_*(Z)$:
\be \sum_{n \geq 0} cl_* (\cF^{\{n\}}) \cdot t^n=\exp \left(- \sum_{r\geq 1}  d_{r*} (\psi_r(cl_*(\cF))\cdot \frac{(-t)^r}{r} \right).
\ee
\ec
In particular, if $Z$ is projective, we recover special cases of the main formulae from \cite{MS}[Cor.1.5], which can now be also derived from (\ref{deg1i}) and (\ref{deg2i}) by specializing the $p_i$'s to $(-1)^{i-1}$. For example, if $cl_*=T_{-y*}$ and $\cF =\bQ^H_Z$, we recover the generating series formula for the degrees $$deg(T_{-y*} ({\bQ^H_Z}^{\{n\}}))=\chi_{-y} ([H^*_c(B(Z,n),\epsilon_n)]),$$
where $B(Z,n) \subset Z^{(n)}$ is the {\it configuration space} of unordered $n$-tuples of distinct points in $Z$, and $\epsilon_n$ is the rank-one local system on $B(Z,n)$ corresponding to a sign representation of $\pi_1(B(Z,n))$ as in \cite{MS}[p.293], compare also with \cite{Go}[Ex.3b] and \cite{Ge}[Cor.5.7].

\medskip

Note also that the specialization $p_1\mapsto 1$ and $p_{i} \mapsto 0$ for all $i \geq 2$, 
corresponds to the evaluation homomorphisms (for all $n \in \bN$) 
$$\frac{1}{n!}ev_{id}=\frac{1}{n!}{\rm Res}^{\Sigma_n}_{id} : {H}_*^{\Sigma_n}(\Zs) \lra {H}_*(\Zs).$$ Then, by forgeting the $\Sigma_n$-action on $\pi_{n*}\cF^{\boxtimes n}$, Corollary \ref{cormaini} specializes to the following result: 
\bc\label{cormain2i}
The following exponential generating series formula holds in the Pontrjagin ring $\mathds{PH}_*(Z)$:
\be  \sum_{n \geq 0} cl_* (\pi_{n*}\cF^{\boxtimes n}) \cdot \frac{t^n}{n!}=\exp \Big( t \cdot cl_*(\cF)  \Big).
\ee
\ec
The above corollary also follows from formula (\ref{m2j}), after a suitable renormalization of the product structure of $\mathds{H}^{\Sigma}_{id,*}(Z)$ in order to make the pushforward $\pi_*:=\oplus \pi_{n,*}:\mathds{H}^{\Sigma}_{id,*}(Z) \to \mathds{PH}_*(Z)$ into a ring homomorphism, see Section \ref{genser} for details.

In particular, if $Z$ is projective, by taking degrees we get exponential generating series formulae for the Euler characteristic and resp. $\chi_y$-polynomial of $H^*(Z^n,\cF^{\boxtimes})$. For example, if $cl_*=T_{-y*}$ and $\cF =\cM$ is a mixed Hodge module on $Z$, we get:
 \be  \sum_{n \geq 0} \chi_{-y} (Z^n,\cM^{\boxtimes n}) \cdot \frac{t^n}{n!}=\exp \Big( t \cdot \chi_{-y}(Z,\cM)  \Big),
\ee
which also follows directly from the K\"unneth formula, e.g., see \cite{MSS}.


\subsection{Twisting by $\Sigma_n$-representations}

Additionally, for a fixed $n$, one can consider the coefficient of $t^n$ in the generating series (\ref{m2i}) for the (equivariant) characteristic classes of all exterior powers $\cF^{\boxtimes n} \in cat^{\Sigma_n}(Z^n)$. Moreover, in this case, one can twist the equivariant coefficients $\cF^{\boxtimes n}$ by a (finite-dimensional) rational $\Sigma_n$-representation $V$, and compute the corresponding equivariant characteristic classes of Lefschetz-type (see Remark \ref{pairing})
\be\label{new1} cl_*(V \otimes \cF^{\boxtimes n};\sigma)=trace_{\sigma}(V) \cdot cl_*(\cF^{\boxtimes n};\sigma),\ee for $\sigma \in \Sigma_n$. By pushing down to the symmetric product $Z^{(n)}$ along the natural map $\pi_n:Z^n \to Z^{(n)}$, we then get by the projection formula the following identity in $$H^{\Sigma_n}_*(Z^{(n)}) \cong H_*(Z^{(n)}) \otimes C(\Sigma_n) \hookrightarrow H_*(Z^{(n)}) \otimes \bQ[p_i, i \geq 1]:$$
\be\label{new2}
cl_*^{\Sigma_n} ( \pi_{n*} (V \otimes \cF^{\boxtimes n})) = \sum_{{\lambda=(k_1,k_2, \cdots) \dashv \ n} } \frac{p_{\lambda}}{z_{\lambda}} \chi_{\lambda}(V)  \cdot \bigodot_{r \geq 1} \left( d_{r*} (\psi_r(cl_*(\cF))) \right)^{k_r}.
\ee
Here, the symbol ${\lambda=(k_1,k_2, \cdots) \dashv \ n}$ denotes the partition  $\lambda=(k_1,k_2, \cdots)$ of $n$ corresponding to a conjugacy class of an element $\sigma\in \Sigma_n$ (i.e., $\sum_r r k_r=n$). We also denote by $z_{\lambda}:=\prod_{r \geq 1} r^{k_r} \cdot k_r!$ the order of the stabilizer of $\sigma$, by $\chi_{\lambda}(V)=trace_{\sigma}(V)$ the corresponding trace, and we set $p_{\lambda}:=\prod_{r \geq 1} p_r^{k_r}$. Finally, $\bigodot$ denotes the Pontrjagin-type product, as defined in Section \ref{genser}. 

If $Z$ is projective, by taking the degree in formula (\ref{new2}), we have the following character formulae generalizing (\ref{deg1i}) and (\ref{deg2i}):
\begin{itemize}
\item[(i)] for $\cF$ a coherent or constructible sheaf, we get 
\be\label{deg1a}  tr_{\Sigma_n} (H^*(Z^n;V \otimes \cF^{\boxtimes n})) =\sum_{{\lambda \dashv \ n} } \frac{p_{\lambda}}{z_{\lambda}} \chi_{\lambda}(V)  \cdot  \chi(H^*(Z;\cF)) ^{\ell(\lambda)},
\ee
with $\chi$ denoting the corresponding Euler characteristic, and for a partition  $\lambda=(k_1,k_2, \cdots)$ of $n$ we let $\ell(\lambda):=k_1+k_2+ \cdots$ be the length of $\lambda$;\newline
\item[(ii)] for $\cM$ a mixed Hodge module on $Z$, we get: 
\be\label{deg2a}  tr_{\Sigma_n} (H^*(Z^n;V \otimes \cM^{\boxtimes n}))=\sum_{{\lambda=(k_1,k_2, \cdots) \dashv \ n} } \frac{p_{\lambda}}{z_{\lambda}} \chi_{\lambda}(V)  \cdot \prod_{r \geq 1} \left( \chi_{-y^r}(H^*(Z;\cM))\right)^{k_r}, 
\ee
with $\chi_y(H^*(Z,\cM))$ the corresponding $\chi_y$-polynomial.
\end{itemize}

Formula (\ref{new2}) is a generalization of Corollary \ref{cormaini}, which one gets back for $V$ the trivial representation. Furthermore, by specializing all the $p_i$'s in equation (\ref{new2}) to the value $1$ (which corresponds to the use of the projectors $(-)^{\Sigma_n}$), one obtains the following identity in $H_*(Z^{(n)})$:
\be\label{new3}
cl_*( (\pi_{n*} (V \otimes \cF^{\boxtimes n}))^{\Sigma_n}) = \sum_{{\lambda=(k_1,k_2, \cdots) \dashv \ n} } \frac{1}{z_{\lambda}} \chi_{\lambda}(V)  \cdot \bigodot_{r \geq 1} \left( d_{r*} (\psi_r(cl_*(\cF)) \right)^{k_r}.
\ee
Note that by letting $V$ be the trivial (resp. sign) representation, formula (\ref{new3}) reduces to Corollary \ref{c1.6} (resp. Corollary \ref{c1.7}). Another important special case of (\ref{new3}) is obtained by choosing $V={\rm Ind}_K^{\Sigma_n}(triv)$, the representation induced from the trivial representation of a subgroup $K$ of $\Sigma_n$, with
$$ (\pi_{n*} (V \otimes \cF^{\boxtimes n}))^{\Sigma_n}
\simeq ({\rm Ind}_K^{\Sigma_n}(triv) \otimes \pi_{n*} (\cF^{\boxtimes n}))^{\Sigma_n}
\simeq (\pi_{n*} (\cF^{\boxtimes n}))^{K} \simeq \pi'_*\left(  (\pi_{*} (\cF^{\boxtimes n}))^{K}  \right).$$
Here $\pi:Z^n \lra Z^n/K$ and $\pi':Z^n/K \to Z^{(n)}$ are the projections factoring $\pi_n$.
In this case, formula (\ref{new3}) calculates the characteristic class $$cl_*((\pi_{n*} (\cF^{\boxtimes n}))^{K})=\pi'_*cl_*((\pi_{*} (\cF^{\boxtimes n}))^{K}).$$ In particular, if $Z$ is projective and we consider the constant Hodge module  $\cF=\bQ^H_Z$, we get at the degree level the following formula for the $\chi_y$-polynomial  of the quotient $Z^n/K$.
\be\label{deg2b}  \chi_{-y} (Z^n/K)=\sum_{{\lambda=(k_1,k_2, \cdots) \dashv \ n} } \frac{1}{z_{\lambda}} \chi_{\lambda}({\rm Ind}_K^{\Sigma_n}(triv))  \cdot \prod_{r \geq 1} \chi_{-y^r}(Z)^{k_r}.
\ee
The corresponding Euler characteristic formula, obtained for $y=1$, is also a special case of Macdonald's formula for the corresponding Poincar\'e polynomial (\cite{Mac}[p.567]).

Finally, by letting $V=V_{\mu} \simeq V^*_{\mu}$ be the (self-dual)  irreducible representation of $\Sigma_n$ corresponding to a partition $\mu$ of $n$, the coefficients $$(\pi_{n*} (V_{\mu} \otimes \cF^{\boxtimes n}))^{\Sigma_n} \simeq 
(V_{\mu} \otimes \pi_{n*} (\cF^{\boxtimes n}))^{\Sigma_n}
=:  S_{\mu}(\pi_{n*} \cF^{\boxtimes n})$$ of the left-hand side of (\ref{new3}) calculate the corresponding Schur functor of $\pi_{n*} \cF^{\boxtimes n} \in cat^{\Sigma_n}(Z^{(n)})$, with
\be\label{new4}
\pi_{n*} \cF^{\boxtimes n} \simeq \sum_{\mu \dashv \ n} V_{\mu} \otimes S_{\mu}(\pi_{n*} \cF^{\boxtimes n}) \in cat^{\Sigma_n}(Z^{(n)}),
\ee
e.g., see Remark \ref{Schur}. These Schur functors generalize the symmetric and alternating powers of $\cF$, which correspond to the trivial and resp. sign representation. 
Note that, by using (\ref{new4}), we get an alternative description of the equivariant classes
$$cl_*^{\Sigma_n}(\pi_{n*} \cF^{\boxtimes n}) \in H^{\Sigma_n}_*(Z^{(n)}) \cong H_*(Z^{(n)}) \otimes C(\Sigma_n) \hookrightarrow H_*(Z^{(n)}) \otimes \bQ[p_i, i \geq 1]$$ in terms of the Schur functions $s_{\mu}:=ch_F(V_{\mu}) \in \bQ[p_i, i \geq 1]$, see \cite{Mc}[Ch.1, Sect.3 and Sect.7]:
\be\label{new5}
cl_*^{\Sigma_n}(\pi_{n*} \cF^{\boxtimes n})=\sum_{\mu \dashv \ n} s_{\mu} \cdot cl_*(S_{\mu}(\pi_{n*} \cF^{\boxtimes n})),
\ee
with $cl_*(S_{\mu}(\pi_{n*} \cF^{\boxtimes n}))$ computed as in (\ref{new3}).

As a concrete example, for $Z$ pure dimensional with coefficients given by the intersection cohomology Hodge module $IC_Z^H$ on $Z$, 
the corresponding Schur functor $S_{\mu}$ of $\pi_{n*} IC^H_{Z^n}$ is given by the twisted intersection cohomology Hodge module $S_{\mu}(\pi_{n*} IC^H_{Z^n})=IC^H_{Z^{(n)}}(V_{\mu})$ with twisted coefficients corresponding to the local system on the configuration space $B(Z,n)$ of unordered $n$-tuples of distinct points in $Z$, induced from $V_{\mu}$ by the group homomorphism $\pi_1(B(Z,n)) \to \Sigma_n$ (compare \cite{MS}[p.293] and \cite{MR}[Prop.3.5]). For $Z$ projective and pure-dimensional, by taking the degrees in (\ref{new3}) for the present choice of coefficients $IC^H_Z$ and representation $V_{\mu}$, we obtain the following identity for the $\chi_y$-polynomial of the twisted intersection cohomology:
\be\label{deg2c}  \chi_{-y} (H^*(Z^{(n)};IC^H_{Z^{(n)}}(V_{\mu}))=\sum_{{\lambda=(k_1,k_2, \cdots) \dashv \ n} } \frac{1}{z_{\lambda}} \chi_{\lambda}(V_{\mu})  \cdot \prod_{r \geq 1} \chi_{-y^r}(H^*(Z;IC^H_Z))^{k_r}.
\ee
Note that results like (\ref{deg2b}) or (\ref{deg2c}) cannot be deduced only from the non-equivariant study of symmetric products as in \cite{CMSSY}. 

\bigskip

We conclude the introduction with a brief discussion of potential applications of our results.

First, the techniques developed here have also been applied by the authors to the study of cohomology representations of external and symmetric products  (see \cite{MS16}), generalizing our previous results from \cite{MS}. 

Secondly, we plan to employ the results of this paper for the study of Hilbert schemes of points on quasi-projective manifolds. In fact, our prior work on symmetric products from \cite{CMSSY} has  already been used for the study of (pushforwards under the Hilbert-Chow morphism of) characteristic classes of Hilbert schemes of points on smooth quasi-projective varieties, see \cite{CMOSY}. But for smooth surfaces, 
results of \cite{CMOSY} may be improved via the McKay correspondence \cite{Kr,Sc1,Sc2}, by using the stronger equivariant results of the present paper. In addition, the present work can also be used for obtaining generating series formulae for the singular Todd classes $td_*(\cF^{[n]})$ of tautological sheaves $\cF^{[n]}$ (associated to a given $\cF \in \co(Z)$) on the Hilbert scheme $Z^{[n]}$ of $n$ points on a smooth quasi-projective algebraic surface $Z$. For degree formulae in this context see, e.g., \cite{WZb}.

One can also use similar techniques for the study of equivariant characteristic classes of the $\Sigma_n$-equivariant Fulton-MacPherson (and other similar) compactifications of configurations spaces of points on a smooth variety $Z$ of any dimension. Here, again, the equivariant results of this paper (and in particular, the twisting by a representation) are needed, while the non-equivariant results from \cite{CMSSY} are not sufficient. For equivariant degree formulae in this context see, e.g., \cite{Ge}.

Finally, techniques and results of this paper can be extended to the context of actions of wreath products $G_n:=G \wr \Sigma_n = G^n \rtimes \Sigma_n$ on external powers $\cF^{\boxtimes n} \in cat^{G_n}(Z^n)$ of a given object $\cF \in cat^G(Z)$ on a $G$-space $Z$, provided that the corresponding key localization formula analogous to (\ref{kl}) holds. This will make the object of future work by the authors. Results of this type for Chern classes already appeared in \cite{Oh2}, while for degree versions see, e.g., \cite{QW,W,WZ,Zh}.

\begin{ack} The authors thank Toru Ohmoto for reading a first version of this paper and for pointing out connections to his work on orbifold Chern classes of symmetric products. The authors are also grateful to the anonymous referees for carefully reading the manuscript, and for their valuable comments and constructive suggestions.

L. Maxim was partially supported by grants from NSF, NSA, by a grant of the Ministry of National Education, CNCS-UEFISCDI project number PN-II-ID-PCE-2012-4-0156, and by a fellowship from the Max-Planck-Institut f\"ur Mathematik,  Bonn.
J. Sch\"urmann was supported by the SFB 878 ``groups, geometry and actions". 
\end{ack}


\section{Delocalized equivariant  theories}\label{dl}

In this section, we introduce the notion of {\it delocalized equivariant theory} of a $G$-space $X$ (with $G$ a finite group), associated to a covariant functor $H$ with compatible cross-product. We also describe the corresponding {\it restriction} and {\it induction} functors, which will play an essential role in the subsequent sections of the paper. 

In the classical context of the usual (co)homology functor, such delocalized theories have been defined in \cite{BBM,BC}, as well as in an unpublished paper of Segal, where they were used for obtaining Riemann-Roch type theorems. For analogues in the context of Deligne-Mumford stacks, see also \cite{Ed, To}. The corresponding orbifold index theorem was developed in \cite{Ka}, by using for the first time the $G$-equivariant cohomology of the inertia space $IX$ (as in (\ref{bc})) of a smooth $G$-space $X$, described in terms of differential forms. The corresponding restriction and induction functors were also studied in this classical context in \cite{QW,Zh}.

\medskip

For simplicity, all spaces in this paper are assumed to be complex quasi-projective, though many constructions in this section apply to other categories of spaces with a finite groups action (e.g., topological spaces or varieties over any base field). 
For such a variety $X$, consider an algebraic action $G\times X \to X$ by a finite group $G$. For any  $g \in G$, we let $X^g$ denote the corresponding  fixed point set.  Let $H$ be a covariant (with respect to isomorphisms) functor to abelian groups, with a compatible cross-product $\boxtimes$ ($\bZ$-linear in each variable), which is commutative,  associative and with a
unit $1_{pt}\in H(pt)$. As main examples used in this paper, we consider the following, with $R$ a commutative ring with unit (e.g., $R=\bZ$, $\bQ$ or $\bC$):
\begin{enumerate}
\item the even degree Borel-Moore homology $H^{BM}_{ev}(X) \otimes R$ of $X$ with coefficients in $R$.
\item Chow groups $CH_*(X) \otimes R$ with $R$-coefficients.
\item Grothendieck group of coherent sheaves $K_0(\co(X)) \otimes R$ with $R$-coefficients.
\end{enumerate}
Other possible choice would be: usual $R$-homology in even degrees $H_{ev}(X) \otimes R$. Since in this section we only need functoriality with respect to isomorphisms, we could also work with cohomological theories, such as the even degree (compactly supported) $R$-cohomology  $H^{ev}_{({c})}(X) \otimes R$ or the Grothendieck group of algebraic vector bundles $K^0(X) \otimes R$ with $R$-coefficients, as used in \cite{QW,Zh}. In this case, the corresponding covariant transformation $g_*$, as used in this paper, is given by $(g^*)^{-1}$, the inverse of the induced pullback under $g$.  If $X$ is smooth, this fits with the folowing Poincar\'e duality isomorphisms:
\be\label{1} 
H_{ev}(X) \otimes R \cong H_c^{ev}(X) \otimes R , \ \ H^{BM}_{ev}(X) \otimes R \cong H^{ev}(X) \otimes R, \ \
K_0({\co(X)}) \otimes R \simeq K^0(X) \otimes R.\ee

\medskip

The disjoint union $\bigsqcup_{g \in G} X^g$ admits an induced $G$-action by $h: X^g \to X^{hgh^{-1}}$, such that the canonical map $$i: \bigsqcup_{g \in G} X^g \to X$$ defined by the inclusions of fixed point sets becomes $G$-equivariant. Therefore, $G$ acts in a natural way on  $\bigoplus_{g \in G}  H(X^g)$.
\bd The {\it delocalized $G$-equivariant theory of $X$ associated to $H$} is the $G$-invariant subgroup of $\bigoplus_{g \in G}  H(X^g)$, namely, \be\label{dh} H^{G}(X):=\left( \bigoplus_{g \in G} \ H(X^g)\right)^G.\ee
This theory is functorial for proper $G$-maps (resp. $G$-equivariant isomorphisms).
 \ed
 
\br  An equivalent interpretation of this delocalized $G$-equivariant theory $H^G(X)$ of a $G$-space $X$ can be obtained by breaking the summation on the right-hand side of (\ref{dh}) into conjugacy classes, i.e.,
\be\label{idex} H^G(X)=\bigoplus_{(g)\in G_*} \left( \bigoplus_{[h]\in G/{Z_G(g)}} h_*\big(H(X^g)^{Z_G(g)}  \big)\right) \cong\bigoplus_{(g)\in G_*} H(X^g)^{Z_G(g)}, \ee
where $G_*$ denotes the set of all conjugacy classes of $G$, and $Z_G(g)$ is the centralizer of
$g\in G$.
\er

\br\label{r2} If $X$ is smooth, then also all fixed-point sets $X^g$ are smooth, so the classical Poincar\'e duality isomorphisms (\ref{1}) induce similar duality isomorphisms $$H_*^G(X) \cong H^*_G(X)$$
between the corresponding delocalized equivariant (co)homology  theories. 
\er
 
 \br If $G$ acts {\it trivially} on $X$ (e.g., $X$ is a point), then 
 \be H^{G}(X) \cong H(X) \otimes C(G),\ee
 where $C(G)$ denotes the free abelian group of $\bZ$-valued  class functions on $G$ (i.e., functions which are constant on the conjugacy classes of $G$).
 \er
 
 \br
 If $G$ is an {\it abelian} group, then:
 \be
 H^G(X)=\bigoplus_{g \in G} H(X^g)^G.
 \ee
 \er
 
 Let us next describe two functors which will be used later.
 
 \bd\label{resdef} ({\it Restriction functor}) \newline
Let $X$ be a $G$-space, as before. For a subgroup $K$ of $G$, the {\it restriction functor} from $G$ to $K$, ${\rm Res}^G_K$, is the group homomorphism $${\rm Res}^G_K:  H^G(X) \to H^K(X)$$
induced by restricting to the $G$-invariant part the projection $$\bigoplus_{g \in G} \ H(X^g) \to \bigoplus_{g \in K} \ H(X^g).$$
Clearly, ${\rm Res}^G_K$ is transitive with respect to subgroups, with ${\rm Res}^G_G$ the identity homomorphism.
 In terms of fixed-point sets of conjugacy classes, i.e., with respect to the isomorphisms:
 $$H^G(X)\cong \bigoplus_{(g)\in G_*} H(X^g)^{Z_G(g)}, \ \ \ \ \  H^K(X)\cong \bigoplus_{(k)\in K_*} H(X^k)^{Z_K(k)},$$
 the restriction factor can be described explicitly as follows (compare with \cite{Zh}[p.4]): if $g \in G$ is not conjugate by elements in $G$ to any element in $K$, then 
 ${\rm Res}^G_K \vert_{H(X^g)^{Z_G(g)}}=0$; otherwise, assume that $g$ is conjugate by elements in $G$ to $k_1,\cdots, k_s \in K$ which have mutually different conjugacy classes in $K$: then $H(X^g)^{Z_G(g)} \cong H(X^{k_i})^{Z_G(k_i)}$ for $i=1,\cdots,s$, and  ${\rm Res}^G_K \vert_{H(X^g)^{Z_G(g)}}$ is given by the direct sum of inclusions $H(X^{k_i})^{Z_G(k_i)} \hookrightarrow H(X^{k_i})^{Z_K(k_i)}$.
 \ed

The following induction functor will be used in the Section \ref{genser} in the definition of Pontrjagin-type products.
\bd\label{inddef} ({\it Induction functor}) \newline
For a $G$-space $X$ as before, and $K$ a subgroup of $G$, the induction from $K$ to $G$, ${\rm Ind}^G_K$, is the group homomorphism (compare with \cite{QW}[p.9]):
\be {\rm Ind}^G_K=\sum_{[g] \in G/K}  g_*(-): H^K(X)\to H^G(X),\ee where the summation is over $K$-cosets of $G$.
In particular, on a $G$-invariant class (i.e., in the image of the restriction functor ${\rm Res}^G_K$) this induction map 
is just multiplication by the index $[G:K]$ of $K$ in $G$. Note that ${\rm Ind}^G_K$ is  transitive for subgroups of $G$, with ${\rm Ind}^G_G$ the identity homomorphism.
In terms of fixed-point sets of conjugacy classes, this induction is given as follows (compare with \cite{Zh}[p.4]): for any conjugacy class $(k)$ in $K$ which intersects the conjugacy class $(g)$ in $G$, we have
\be\label{one} {\rm Ind}^G_K=  \sum_{[h]\in Z_G(k)/{Z_K(k)}} h_*(-) :
H(X^k)^{Z_K(k)}  \to H(X^k)^{Z_G(k)} \cong H(X^g)^{Z_G(g)},\ee
so on a $G$-invariant class this is just  multiplication
by the index $[Z_G(k):Z_K(k)]$.
\ed

\br\label{expl} In terms of the above induction functors, the identification (\ref{idex}) is given by:
$$\bigoplus_{(g)\in G_*} {\rm Ind}_{Z_G(g)}^G: \bigoplus_{(g)\in G_*} H(X^g)^{Z_G(g)} \to H^G(X), $$
where $ {\rm Ind}_{Z_G(g)}^G: H(X^g)^{Z_G(g)} \to H^G(X)$ is the restriction of ${\rm Ind}_{Z_G(g)}^G$ to the direct summand 
$$ H(X^g)^{Z_G(g)} \hookrightarrow H^{Z_G(g)}(X)$$ coming from the $Z_G(g)$-equivariant direct summand $H(X^g) \subset \bigoplus_{h \in Z_G(g)} H_*(X^h)$.
\er

\br Over a point space, the above functors reduce in many cases to the classical restriction and induction functors from the representation theory of finite groups.
\er


\subsection{Compatibilities with cross-product}\label{1.1}
Assume $G$ acts on $X$, with $g\in G$ and $K\subset G$ a subgroup, and similarly for $G'$ acting on $X'$, with $g'\in G'$ and $K' \subset G'$ a subgroup.
Then $(X\times X')^{g\times g'}=X^g \times X'^{g'}$,
$Z_{G\times G'}(g\times g')=Z_G(g)\times Z_{G'}(g')$,
as well as $G\times G'/K\times K'= G/K \times G'/K'$,
and similarly for the quotient of centralizers as above.

Then all products $\boxtimes: H(X^g)\times H(X'^{g'})\to H(X^g\times X'^{g'})$ 
induce by the functoriality of $\boxtimes$ a corresponding
commutative and associative cross-product 
$$\boxtimes: H^G(X)\times H^{G'}(X') \to H^{G\times G'}(X\times X')$$
(with unit $1_{pt} \in H^{\{id\}}(pt)$, for  $\{id\}$ denoting the trivial group).
Moreover, this product is compatible with the restriction and induction functors, i.e.,
\be {\rm Ind}^{G\times G'}_{K\times K'}( - \boxtimes -) = 
{\rm Ind}^G_K(-) \boxtimes {\rm Ind}^{G'}_{K'}(-)\ee
and
\be {\rm Res}^{G\times G'}_{K\times K'}( - \boxtimes -) = 
{\rm Res}^G_K(-) \boxtimes {\rm Res}^{G'}_{K'}(-).\ee

Finally, the above facts about cross-product and restriction functors can be used to define a pairing:
$$C(G) \times H^G(X) \overset{\cdot}{\lra} H^G(X)$$
by $$H^G(pt) \times H^G(X) \overset{\boxtimes}{\lra} H^{G\times G}(pt \times X) \overset{{\rm Res}}{\lra} H^G(X),$$
with $pt \times X \cong X$, and ${\rm Res}$ denoting the restriction functor for the diagonal subgroup $G \hookrightarrow G \times G$.

\br\label{Rid} The distinguished unit element $id \in G$ 
gives the direct summand 
$$H^G_{id}(X):=H(X)^G \subset H^G(X),$$
i.e., the $G$-invariant subgroup $H(X)^G$ of $H(X)$. 
This direct  summand is compatible with restriction, induction and
induced cross-products. If the functor $H$ is also covariantly functorial for closed embeddings, we get a pushforward
for the closed fixed point inclusions $i_g: X^g \hookrightarrow X$, i.e., 
$$i_{g*}: H(X^g)\to H(X),$$ and a group homomorohism 
$$sum_G:=\sum_g i_{g*}: H^G(X)\to H^G_{id}(X)=H(X)^G \subset H(X).$$
Note that this homomorphism commutes with induction and cross-products.
\er

\br\label{r3}
If $X$ is smooth, the induction-restriction functors, as well as their compatibilities with cross-products are also compatible with Poincar\'e duality for (co)homology as in Remark \ref{r2}.
\er


\section{Generating series for symmetric group actions on external products}\label{genserg}

In this section, we describe a very general generating series formula for symmetric group actions on external products, which should be regarded as a far-reaching generalization of a well-known identity of symmetric functions. In section \ref{exabs}, we give applications of this abstract generating series formula in the context of orbifold cohomology and resp. localized $K$-theory.

\medskip

Let $Z$ be a quasi-projective variety, with the symmetric group $\Sigma_n$ acting on  the cartesian product $Z^n$ of $n \geq 0$ copies of $Z$ by the natural permutation action. For our generating series formula, it is important to look at all groups $H^{\Sigma_n}(Z^n)$ simultaneously.
Let 
\be
\mathds{H}^{\Sigma}(Z):=\bigoplus_{n \geq 0} H^{\Sigma_n}(Z^n) \cdot t^n
\ee
be the commutative graded $\bZ$-algebra (with unit) with product $$\odot:={\rm Ind}^{\Sigma_{n+m}}_{\Sigma_n \times \Sigma_m}(\cdot \boxtimes \cdot)$$ induced from the external product by induction. Here, $\bigoplus_{n \geq 0}H^{\Sigma_n}(Z^n)$ becomes a commutative graded ring with product $\odot$, and we view the completion $\mathds{H}^{\Sigma}(Z)$ as a subring of the formal power series ring $\bigoplus_{n \geq 0} H^{\Sigma_n}(Z^n)[[t]]$. 

The algebra $\mathds{H}^{\Sigma}(Z)$ is, in addition, endowed with  {\it creation operators} $$\mathfrak{a}_r: H(Z) \to H^{\Sigma_r}(Z^r),$$ $r \geq 1$, which allow us to transport elements from $H(Z)$ to the delocalized groups $H^{\Sigma_r}(Z^r)$. These are defined as follows:
 if $\sigma_r=({r})$ is an $r$-cycle in $\Sigma_r$, then $\mathfrak{a}_r$ is the composition
$$\mathfrak{a}_r:H(Z) \overset{\cdot r}{\to} H(Z) \cong H((Z^r)^{\sigma_r})^{Z_{\Sigma_r}(\sigma_r)} \hookrightarrow H^{\Sigma_r}(Z^r),$$
where $\langle \sigma_r \rangle =Z_{\Sigma_r}(\sigma_r)$ acts trivially on $(Z^r)^{\sigma_r}$, and therefore  also on $H((Z^r)^{\sigma_r})$. The role of multiplication by $r$ in the definition of creation operator will become clear later on, e.g., in the proof of Theorem \ref{main2} below.
The creation operator $\mathfrak{a}_r$ can be re-written as $$\mathfrak{a}_r:=r\cdot {\rm Ind}^{\Sigma_r}_{\langle \sigma_r \rangle}\circ i_r,$$ with $$i_r: H(Z)\simeq H((Z^r)^{\sigma_r})^{\langle \sigma_r \rangle } \subset H^{<\sigma_r>}(Z^r).$$ Here, the last inclusion is just a direct summand, because $\langle \sigma_r \rangle$ is abelian. In the following we omit to mention $i_r$ explicitly.

\medskip

Let $\sigma \in \Sigma_n$ have cycle partition $\lambda=(k_1, k_2, \cdots )$, i.e.,  $k_r$ is the number of length $r$ cycles in $\sigma$ and $n=\sum_r r \cdot k_r$. 
Then 
$$(Z^n)^{\sigma} \simeq  \prod_r \left( (Z^r)^{\sigma_r} \right)^{k_r} \simeq \prod_r \Delta_r(Z)^{k_r} \simeq Z^{k_1+k_2+\cdots} \ ,$$ 
where $\sigma_r$ denotes as above a cycle of length $r$ in $\Sigma_n$, and $\Delta_r(Z)$ is the diagonal in $Z^r$, i.e., the image of the diagonal map $\Delta_r:Z \to Z^r$. 

Let us now choose a sequence $\underline{b}=(b_1, b_2,\cdots)$ of elements $b_r\in H(Z)$, $r\geq 1$, and associate to a conjugacy class represented by  $\sigma \in \Sigma_n$ of type $(k_1,k_2,\cdots)$ the element ${\underline{b}}^{(\sigma)} \in H^{\Sigma_n}(X^n)$ corresponding to
$$\boxtimes_r (b_r)^{\boxtimes k_r} \in H(\prod_r Z^{k_r})\simeq H((Z^n)^{\sigma}),$$ as it will be explained below.
Recall that $Z_{\Sigma_n}(\sigma)$ is a product over $r$ of semidirect products of $\Sigma_{k_r}$ with $\langle \sigma_r \rangle^{k_r}$, that is, \be Z_{\Sigma_n}(\sigma) \cong \prod_r \Sigma_{k_r} \ltimes \bZ_r^{k_r}\ee
(with $\sigma_r$ denoting as before an $r$-cycle). The group $\bZ_r^{k_r} \cong  \langle \sigma_r \rangle^{k_r}$ acts trivially on $Z^{k_r}$, whereas $\Sigma_{k_r}$ permutes the corresponding $Z$-factors of $Z^{k_r}$ (compare \cite{Zh}[p.8]). By commutativity and associativity of the cross-product $\boxtimes$, it follows that $\boxtimes_r (b_r)^{\boxtimes k_r}$ is invariant under $Z_{\Sigma_n}(\sigma)$, so it indeed defines an element \be\label{elem}{\underline{b}}^{(\sigma)}={\rm Ind}^{\Sigma_n}_{Z_{\Sigma_n}(\sigma)}\left( \boxtimes_r (b_r)^{\boxtimes k_r} \right) \in  H^{\Sigma_n}(Z^n),\ee
with induction defined as in Remark \ref{expl}.
Moreover, for $\sigma \in \Sigma_n$ and $\sigma' \in \Sigma_{m}$, we have:
\be {\underline{b}}^{(\sigma)} \odot {\underline{b}}^{(\sigma')}=
{\underline{b}}^{(\sigma \times \sigma')} \in H^{\Sigma_{n+m}}(Z^{n+m}).
\ee

\medskip

In what follows, we assume that the functor $H$ takes values in $R$-modules, with $R$ a commutative $\bQ$-algebra (otherwise, work with $\mathds{H}^{\Sigma}(Z) \otimes R$). It follows that $\mathds{H}^{\Sigma}(Z)$ is also a commutative graded  $\bQ$-algebra.
Note that one can also switch between covariant and contravariant notions, e.g. between homology and cohomology
by Poincar\'e duality, if $X$ is smooth.

The main result of this section is the following generating series formula:
\bt\label{main2} 
With the above notations, the following generating series formula holds in the $\bQ$-algebra  $\mathds{H}^{\Sigma}(Z)$:
\be\label{m2} \sum_{n \geq 0} \left( \sum_{(\sigma) \in (\Sigma_n)_*} \underline{b}^{(\sigma)} \right) \cdot t^n=\exp \left( \sum_{r\geq 1} \mathfrak{a}_r (b_r)\cdot \frac{t^r}{r} \right),
\ee
where $(\Sigma_n)_*$ denotes the set of conjugacy classes of $\Sigma_n$. \et

\begin{proof}
We have the following string of equalities in the $\bQ$-algebra $(\mathds{H}^{\Sigma}(Z),\odot)$:
\begin{equation}
\begin{split}
\exp\left( \sum_{r=1}^{\infty} x_r \frac{t^r}{r} \right)  &=
\prod_{r=1}^{\infty} \exp \left(  x_r \frac{t^r}{r} \right) \\
&=\prod_{r=1}^{\infty} \sum_{k_r=0}^{\infty} \left(  x_r \frac{t^r}{r} \right)^{k_r} \frac{1}{k_r !} \\
&=\sum_{N \geq 0} \sum_{k_1,\cdots,k_N} \frac{x_1^{k_1} \cdots x_N^{k_N}}{k_1!\cdots k_N!} \prod_{r=1}^N \left( \frac{t^r}{r} \right)^{k_r} \\
&=\sum_{N \geq 0} \sum_{k_1,\cdots,k_N} \frac{x_1^{k_1} \cdots x_N^{k_N}}{k_1!\cdots k_N!} \frac{t^{k_1+2k_2+\cdots +Nk_N}}{1^{k_1}\cdots N^{k_N}}\\
&=\sum_{m=0}^{\infty} t^m \sum_{k_1+2k_2+\cdots + Nk_N=m} \frac{x_1^{k_1} \cdots x_N^{k_N}}{k_1!\cdots k_N! 1^{k_1}\cdots N^{k_N}} \\
&=\sum_{m=0}^{\infty} t^m \sum_{k_1+2k_2+\cdots + Nk_N=m} 
\prod_{r=1}^N  \frac{x_r^{k_r}}{k_r! r^{k_r}}  \\
\end{split}
\end{equation}
Note that the sum over $k_1+2k_2 + \cdots Nk_N=m$ corresponds to a summation over the cycle classes $(\sigma)$ in $\Sigma_m$ given by $\prod_r {\sigma_r}^{k_r}$, for $\sigma_r=({r})$ an $r$-cycle in $\Sigma_r$. 
In our case, we take:
$$x_r=\mathfrak{a}_r(b_r)={\rm Ind}^{\Sigma_r}_{\langle \sigma_r \rangle}(r\cdot b_r).$$ 

All products (and powers) above are with respect to the multiplication $\odot$ in $\mathds{H}^{\Sigma}(Z)$, which is  defined via cross-product and induction. In particular, 
 $$\prod_{r=1}^N (x_r)^{k_r}/r^{k_r}={\rm Ind}^{\Sigma_m}_{\prod_r (\Sigma_r)^{k_r}}  (\boxtimes_{r=1}^N (x_r)^{\boxtimes k_r}/r^{k_r} ).$$
Moreover, by using the compatibility of induction with $\boxtimes$, $\bZ$-linearity and transitivity, we have: $${\rm Ind}^{\Sigma_m}_{\prod_r (\Sigma_r)^{k_r}}  (\boxtimes_{r=1}^N (x_r)^{\boxtimes k_r}/r^{k_r} )= {\rm Ind}^{\Sigma_m}_{\prod_r \langle \sigma_r \rangle^{k_r}} ( \boxtimes_r (b_r)^{\boxtimes k_r} )= {\rm Ind}^{\Sigma_m}_{Z_{\Sigma_m}(\sigma)}\circ  {\rm Ind}^{Z_{\Sigma_m}(\sigma)}_{\prod_r \langle \sigma_r \rangle^{k_r}} ( \boxtimes_r (b_r)^{\boxtimes k_r}),$$ where, as before, $\sigma$ is a representative of the cycle type $(k_1,k_2,\cdots)$. But, as already mentioned, $Z_{\Sigma_m(\sigma)}$ acts trivially on $ \boxtimes_r (b_r)^{\boxtimes k_r}$, so ${\rm Ind}^{Z_{\Sigma_m(\sigma)}}_{\prod_r \langle \sigma_r \rangle^{k_r}}$ is just multiplication by the index $[Z_{\Sigma_m(\sigma)}:\prod_r \langle
 \sigma_r \rangle^{k_r}]=\prod_r k_r!$.

Altogether, we get
$$\prod_{r=1}^N (x_r)^{k_r}/k_r! r^{k_r}= {\underline b}^{(\sigma)},$$
which finishes the proof.
\end{proof}


\subsection{Examples}\label{exabs} Let us now explain some special cases of Theorem \ref{main2} in the cohomological language, which in some situations are already available in the literature. Our main applications, to equivariant characteristic classes for singular spaces, will be given later on, in Section \ref{gensersp}, after we develop the necessary background. 

\subsubsection{Orbifold cohomology}
Here we work with $H(X):=H^{ev}(X) \otimes \bQ$, the (even degree) rational cohomology functor. For $X$ smooth, our notion of $H^G(X)$ corresponds to the even degree orbifold cohomology $H^{2*}_{\rm orb}(X/G)$, as used for example in \cite{QW}.

For $Z$ a quasi-projective complex variety, and for a given $\gamma \in H(Z)$, let $b_r:=\gamma$, for all $r \geq 1$. Then ${\underline b}^{(\sigma)}$ corresponds to $\gamma^{\boxtimes \ell(\sigma)} \in H((Z^n)^{\sigma})$ for $\sigma \in \Sigma_n$  of cycle type $(k_1, k_2, \cdots)$, and $\ell(\sigma):=\sum_r k_r$ the length of the partition associated to $\sigma$. Following \cite{QW}, we set:
\be \eta_n(\gamma):= \sum_{(\sigma) \in (\Sigma_n)_*} 
{\rm Ind}^{\Sigma_n}_{Z_{\Sigma_n(\sigma)}}\left( \gamma^{\boxtimes \ell(\sigma)} \right)
=  \sum_{(\sigma) \in (\Sigma_n)_*} \underline{b}^{(\sigma)}.
\ee
Then our formula (\ref{m2}) specializes to the following result:
\be\label{m2c1} \sum_{n \geq 0} \eta_n(\gamma) \cdot t^n=\exp \left( \sum_{r\geq 1} \mathfrak{a}_r (\gamma)\cdot \frac{t^r}{r} \right).
\ee
For $X$ smooth, this fits with the formula stated after Definition 3.2 in \cite{QW}. However, our proof is purely formal, so it applies to any topological space, as well as to algebraic varieties with rational Chow groups for $H$.


\subsubsection{Localized $K$-theory}\label{loc1}
Here we work with $H(X):=K^0(X) \otimes \bC$, the complexified Grothendieck group of algebraic vector bundles on $X$. For a finite group $G$ acting algebraically on a quasi-projective complex variety $X$, we define a localization map of Lefschetz-type
$$L^G:K_G^0(X) \lra H^G(X)$$ on the Grothendieck group of $G$-equivariant algebraic vector bundles as a direct sum of transformations:
$$L(g):[W] \in K^0_{G} (X)
\mapsto [W|_{X^g}] \in K^0_{\langle g \rangle}(X^g)\simeq 
 K^0 (X^g)\otimes {\it Rep}_{\bC}(\langle g \rangle) \mapsto K^0 (X^g) \otimes \bC,$$
 where the last map is induced by taking the trace against $g\in G$. Here, ${\it Rep}_{\bC}(\langle g \rangle)$ is the Grothendieck group of complex representations of the group $\langle g \rangle$, and the isomorphism in the above definition holds since $\langle g \rangle$ acts trivially on the fixed-point $X^g$ (e.g., this fact follows from \cite{D}[(1.3.4)], compare also \cite{Ge}). Note that, by construction, $L^G$ commutes with cross-products as in Section \ref{1.1}
 
 \medskip
 
For $Z$ a quasi-projective complex variety and an algebraic vector bundle $V$ on $Z$, we get the $\Sigma_n$-equivariant vector bundle $V^{\boxtimes n}$ on $Z^n$. Let $\sigma \in \Sigma_n$ be of cycle type $(k_1,k_2,\cdots)$. Then, by the multiplicativity of $L(\sigma)$, we get:
$$L(\sigma)=\boxtimes_r L(\sigma_r)^{\boxtimes k_r}.$$ So it suffices to understand the transformations $L(\sigma_r)$, for all $r$-cycles ($r \geq 1$).

For $\sigma_n=(n)$ an $n$-cycle, we have that 
\be L(\sigma_n)[V^{\boxtimes n}]=\psi_n(V) \in K^0(Z) \otimes \bQ
\ee
under the identification $(Z^n)^{\sigma_n}=\Delta_n(Z) \simeq Z$, with $$\Delta_n^*(V^{\boxtimes n})=V^{\otimes n}.$$ Here $\psi_n$ denotes the $n$-th Adams operation defined by Atiyah in the topological context \cite{At} and e.g., Nori \cite{N}[Lem.3.2] in the algebraic geometric context. Note that we can work here with rational coefficients, since characters of symmetric groups are integer-valued.

If we choose $b_r:=\psi_r(V)$, for all $r \geq 1$, our main formula  (\ref{m2}) specializes to the following generating series identity:
\be\label{m2c2} \sum_{n \geq 0}  L^{\Sigma_n}(V^{\boxtimes n} ) \cdot t^n=\exp \left( \sum_{r\geq 1} \mathfrak{a}_r (\psi_r(V))\cdot \frac{t^r}{r} \right),
\ee
since, by {\it multiplicativity} and {\it conjugacy invariance} of $L(\sigma)$, we have that $$L^{\Sigma_n}(V^{\boxtimes n})=\sum_{(\sigma) \in (\Sigma_n)_*} {\rm Ind}^{\Sigma_n}_{Z_{\Sigma_n}(\sigma)}\left( \boxtimes_r L(\sigma_r)^{\boxtimes k_r} \right)= \sum_{(\sigma) \in (\Sigma_n)_*} \underline{b}^{(\sigma)}.$$
The same proof applies in the topological context, for topological $K$-theory, in which case we obtain a special situation (for $G$ the identity group) of \cite{W}[Prop.4]. Note that \cite{W} uses the identification $L^G \otimes \bC :K^0_G(X) \otimes \bC \simeq H^G(X)$.

Similarly, one can work with algebraic varieties over any base field of characteristic zero, with $L(g)$ the corresponding Lefschetz transformation of Baum-Fulton-Quart \cite{BFQ}.


\subsubsection{Localized Grothendieck groups of constructible sheaves}\label{loc2}
Here we work with $$H(X):=K_0({\cons(X)}) \otimes \bC,$$ the complexified Grothendieck group of (algebraically) constructible sheaves of complex vector spaces. For a finite group $G$ acting algebraically on a quasi-projective complex variety $X$, we define a localization map of Lefschetz-type
$$L^G:K^G_0({\cons(X)}) \lra H^G(X)$$ on the Grothendieck group of $G$-equivariant (algebraically) constructible sheaves of complex vector spaces as a direct sum of transformations:
\begin{align*}L(g):[\cF] \in K_0^{G} ({\cons(X)})
&\mapsto [\cF|_{X^g}]  \in K_0^{\langle g \rangle}({\cons(X^g)}) \\
& \simeq 
 K_0 ({\cons(X^g)})\otimes {\it Rep}_{\bC}(\langle g \rangle) \mapsto K_0 ({\cons(X^g)}) \otimes \bC,\end{align*}
 where the last map is induced, as before,  by taking the trace against $g\in G$. The isomorphism in the above definition holds since $\langle g \rangle$ acts trivially on the fixed-point $X^g$ and ${\cons(X^g)}$ is an abelian $\bC$-linear category (e.g., see as before \cite{D}[(1.3.4)], compare also \cite{Ge}). Note that, by construction, $L^G$ commutes as before with cross-products as in Section \ref{1.1}
 
 \medskip
 
For $Z$ a quasi-projective complex variety and an (algebraically) constructible sheaf $\cF$ on $Z$, we get the $\Sigma_n$-equivariant  (algebraically) constructible sheaf $\cF^{\boxtimes n}$ on $Z^n$. For $\sigma_n=(n)$ an $n$-cycle, we have that 
\be L(\sigma_n)[\cF^{\boxtimes n}]=\psi_n(V) \in K_0({\cons(Z)}) \otimes \bQ
\ee
under the identification  $(Z^n)^{\sigma_n}=\Delta_n(Z) \simeq Z$, with $$\Delta_n^*(\cF^{\boxtimes n})=\cF^{\otimes n}.$$ Here $\psi_n$ denotes the $n$-th Adams operation of the pre-lambda ring structure on $K_0 ({\cons(Z)})$ induced from the symmetric monoidal tensor product $\otimes$ of constructible sheaves, as in \cite{MS}[Lem.2.1].
As before, we can work here with rational coefficients.

If we choose $b_r:=\psi_r(\cF)$, for all $r \geq 1$, our main formula  (\ref{m2}) specializes as above to the following generating series identity:
\be\label{m2c3} \sum_{n \geq 0}  L^{\Sigma_n}(\cF^{\boxtimes n} ) \cdot t^n=\exp \left( \sum_{r\geq 1} \mathfrak{a}_r (\psi_r(\cF))\cdot \frac{t^r}{r} \right).
\ee

\subsubsection{Frobenius character}\label{Fr}
Specializing to a point space $X$, the above localized theories for vector bundles and resp. constructible sheaves reduce to the classical character theory of a finite group $G$:
$$tr_G:{\it Rep}_{\bC}(G) \lra C(G) \otimes \bC$$
$$[V] \mapsto \{trace_g(V), g \in G\},$$
with ${\it Rep}_{\bC}(G)$ the Grothendieck group of complex representations of $G$, and $C(G)$ the free abelian group of $\bZ$-valued class functions on $G$. The trace $trace_g$ is of course multiplicative and conjugacy invariant. 

For symmetric groups, we can work again with rational coefficients, and get an algebra homomorphism 
$$tr_{\Sigma}:{\it Rep}_{\bC}(\Sigma):=\bigoplus_n {\it Rep}_{\bC}(\Sigma_n) \lra C(\Sigma)\otimes \bQ:=\bigoplus_n C(\Sigma_n) \otimes \bQ$$
with respect to the classical induction product: $\odot:={\rm Ind}^{\Sigma_{n+m}}_{\Sigma_n \times \Sigma_m}(\cdot \boxtimes \cdot)$ for representations and resp. characters, see e.g., \cite{Mc}[Ch.I,Sect.7]. This homomorphism can be composed with the {\it Frobenius character} $$ch_F: C(\Sigma)\otimes \bQ:=\bigoplus_n C(\Sigma_n) \otimes \bQ \overset{\simeq}{\lra} \bQ[p_i, i \geq 1]=:\Lambda \otimes \bQ$$
to the graded ring of $\bQ$-valued symmetric functions in infinitely many variables $x_m$ ($m \in \bN$), with $p_i:=\sum_m x_m^i$ the $i$-th power sum function. On $C(\Sigma_n) \otimes \bQ$, $ch_F$ is defined by:
\be\label{Frob}
ch_F (f):=\frac{1}{n!} \sum_{\sigma \in \Sigma_n} f(\sigma) \psi(\sigma),
\ee
with $$\psi(\sigma):=\prod_r p_r^{k_r}$$ for $\sigma$  of cycle type $(k_1,k_2, \cdots)$, e.g., see \cite{Mc}[Ch.I,Sect.7]. For example, if $f$ is the indicator function of the conjugacy class of the $n$-cycle $\sigma_n$  in $\Sigma_n$, then $ch_F (f)=\frac{1}{n} p_n$ since $n=|Z_{\Sigma_n}(\sigma_n)|$.
In particular, the creation operator $\mathfrak{a}_r :\bQ \to \bQ[p_i, i \geq 1]=\Lambda \otimes \bQ$ is (up to the Frobenius isomorphism) given by multiplication with $p_r$, which also motivates the use of multiplication by $r$ in the definition of our creation operator in Section \ref{genserg}.

If we choose $b_r:=1 \in \bQ$, for all $r \geq 1$, our main formula  (\ref{m2}) specializes to the well-known identity of symmetric functions (e.g., see the proof of \cite{Mc}[(2.14)]):
\be\label{m2c4} H(t):=\sum_{n \geq 0}  h_n t^n=\exp \left( \sum_{r\geq 1}  p_r \cdot \frac{t^r}{r} \right),
\ee
with $h_n=ch_F(1_{\Sigma_n})$ the $n$-th complete symmetric function (see \cite{Mc}[p.113]), and $1_{\Sigma_n}:=tr_{\Sigma_n}(triv_n)$ the identity character of the trivial representation $triv_n$ of $\Sigma_n$.


\section{(Equivariant) Pontrjagin rings for symmetric products}\label{genser}

Let $Z$ be a quasi-projective variety, and denote by $\Zs$ its $n$-th symmetric product, i.e., the quotient of the product $Z^n$ of $n$ copies of $Z$ by the natural action of the symmetric group on $n$ elements, $\Sigma_n$. Let $\pi_n:Z^n \to \Zs$ denote the natural projection map. 

In this section, the functor $H$ from Section \ref{dl} is required in addition to be covariant (at least) for {\it finite maps} such as $\pi_n$ or the closed embedding $i_{\sigma}:(Z^n)^{\sigma} \hookrightarrow Z^n$, for $\sigma \in \Sigma_n$. We will carry over the assumption that $H$ takes values in $R$-modules, with $R$ a commutative $\bQ$-algebra.

Besides $\mathds{H}^{\Sigma}(Z):=\bigoplus_{n \geq 0} H^{\Sigma_n}(Z^n) \cdot t^n$, here we consider other structures of commutative graded $\bQ$-algebra with units, defined in terms of symmetric and respectively external products of $Z$:
\bn
\item[(a)] On $$\mathds{PH}(Z):=\bigoplus_{n \geq 0} H(\Zs) \cdot t^n=\prod_{n \geq 0} H(\Zs)$$ there is the Pontrjagin ring structure, with multiplication $\odot$ induced from the maps $$\Zs \times Z^{(m)} \to Z^{(m+n)},$$ see \cite{CMSSY}[Definition 1.1] for more details. Here, $\bigoplus_{n \geq 0} H(\Zs) $ becomes a commutative graded ring with product $\odot$, and we view the completion $\mathds{PH}(Z)$ as a subring of the formal power series ring $\bigoplus_{n \geq 0} H(\Zs)[[t]]$.
\item[(b)] On 
\begin{align*} \mathds{PH}^{\Sigma}(Z):=\bigoplus_{n \geq 0} H^{\Sigma_n}(\Zs) \cdot t^n &\cong \bigoplus_{n \geq 0} \left( H(\Zs) \otimes C(\Sigma_n) \otimes \bQ \right) \cdot t^n \\&\hookrightarrow \mathds{PH}(Z) \otimes \left( C(\Sigma) \otimes \bQ \right) \end{align*}
there is a product induced from that of the Pontrjagin product in the $H$-factor and the induction product for class functions. Via the Frobenius character identification $$ch_F:C(\Sigma) \otimes \bQ \simeq \bQ[p_i \ , i \geq 1], $$ we can also view $\mathds{PH}^{\Sigma}(Z)$ as a graded subalgebra of $\mathds{PH}(Z) \otimes \bQ[p_i \ , i \geq 1]$, and with $i$-th power sum $p_i$ regarded as a degree $i$ variable. 
\item[({c})] By Remark \ref{Rid}, the direct summand $$\mathds{H}^{\Sigma}_{id}(Z):=\bigoplus_{n \geq 0} H_{id}^{\Sigma_n}(Z^n) \cdot t^n \subset \mathds{H}^{\Sigma}(Z)$$ corresponding to the identity component is a subring, so that 
$$\bigoplus_{n \geq 0} sum_{\Sigma_n}:\mathds{H}^{\Sigma}(Z) \lra \mathds{H}^{\Sigma}_{id}(Z)$$
is a ring homomorphism. With respect to the Frobenius homomorphism, it is more natural to use the averaging homomorphisms $av_n:=\frac{1}{n!} sum_{\Sigma_n}:H^{\Sigma_n}(Z^n) \to H^{\Sigma_n}_{id}(Z^n)$. Then the graded group homomorphism $$av:=\bigoplus av_n: \mathds{H}^{\Sigma}(Z)\to \mathds{H}^{\Sigma}_{id}(Z)$$ becomes a graded algebra homomorphism if we introduce on
$\mathds{H}^{\Sigma}_{id}(Z)$ the {\it twisted product}
$$\widetilde{\odot}:= \frac{n!m!}{(n+m)!} \odot :
H^{\Sigma_n}_{id}(Z^n)\times H^{\Sigma_m}_{id}(Z^m) \to
H^{\Sigma_{n+m}}_{id}(Z^{n+m}).$$
With this twisted product, we also have a Frobenius-type ring homomorphism 
$$av_F: \mathds{H}^{\Sigma}(Z)\to \mathds{H}^{\Sigma}_{id}(Z)\otimes \bQ[p_i, i\geq 1]$$
given by
$$\frac{1}{n!} \sum_{\sigma \in \Sigma_n} i_{\sigma *}\cdot \psi(\sigma):
H^{\Sigma_n}(Z^n)\to H^{\Sigma_n}_{id}(Z^n)\otimes \bQ[p_i, i\geq 1],$$
with $\psi(\sigma)$ as in the Frobenius homomorphism $ch_F$ of equation (\ref{Frob}).
\en
These structures are related by homomorphisms of commutative graded $\bQ$-algebras, fitting into 
the following commutative diagram:
\be\label{two}
\xymatrix{
\mathds{H}^{\Sigma}(Z):=\bigoplus_{n \geq 0} H^{\Sigma_n}(Z^n) \cdot t^n \ar[d]^{\pi_*=\oplus_n \pi_{n*}}  \ar[r]^{\ \ \ \ \ \ av_F} &  \mathds{H}^{\Sigma}_{id}(Z)\otimes \bQ[p_i, i\geq 1]  \ar[d]^{\pi_* \otimes id} \\
 \mathds{PH}^{\Sigma}(Z):=\bigoplus_{n \geq 0} H^{\Sigma_n}(\Zs) \cdot t^n \ar[d]^{\oplus_n \frac{1}{\vert \Sigma_n \vert} \sum_{\sigma} ev_{\sigma}}  \ar[r]  &    \mathds{PH}(Z) \otimes \bQ[p_i , \ i \geq 1] \ar[d]^{p_i=1} \\ 
 \mathds{PH}(Z):=\bigoplus_{n \geq 0} H(\Zs) \cdot t^n \ar[r]^{\ \ \ \ \ \ \ \ \ \ \ \ \ =} &   \mathds{PH}(Z)  } 
\ee
with $ev_{\sigma}:C(\Sigma_n) \otimes \bQ \to \bQ$ the evaluation map at $\sigma \in \Sigma_n$.

\medskip

Let $d_r:=\pi_r \circ \Delta_r$ be the composition $d_r:Z \to Z^r \to Z^{({r})}$ of the natural projection $\pi_r:Z^r \to Z^{({r})}$ with the diagonal embedding $\Delta_r:Z \to Z^r$. Then the creation operator $\mathfrak{a}_r$ satisfies the identities:
\be\label{cr} \pi_{r*} \circ \mathfrak{a}_r=p_r \cdot d_{r*}, \ \ \ av_F \circ \mathfrak{a}_r=p_r \cdot \Delta_{r*}, \ \ \ {\rm and} \ \ \ av \circ \mathfrak{a}_r=\Delta_{r*}.\ee
This generalizes the corresponding relation between $\mathfrak{a}_r$ and $p_r$ discussed at the end of Section \ref{Fr}.

\br
Under the assumptions from the beginning of this section, the commutative diagram (\ref{two}) is  functorial in $Z$ for finite maps. If, moreover, the functor $H$ and the cross-product are functorial for proper morphisms, then  (\ref{two})  is also functorial for such morphisms. In particular, for $Z$ compact, we can push down our generating series formulae (such as (\ref{m2}))
to a point to obtain (equivariant) degree formulae. Finally, the diagram (\ref{two}) is compatible with natural transformations of such functors.
\er

\subsection{Example: Constructible functions and Orbifold Chern classes}
In this example, we explain how our main result of Theorem \ref{main2} can be used to reprove Ohmoto's generating series identities for canonical contructible functions (\cite{Oh2}[Prop.3.9]) and orbifold Chern classes (\cite{Oh2}[Thm.1.1]). These are generalized class versions of the celebrated Hirzebruch-H\"ofer (or Atiyah-Segal) formula for the orbifold Euler characteristic of symmetric products. Recall that for a (compact) $G$-space $X$ as before, the orbifold Euler characteristic $\chi(X,G)$ is defined as:
\be\label{oe} \chi(X,G)=\frac{1}{|G|} \sum_{gh=hg} \chi(X^g \cap X^h)=\sum_{(g) \in G_*} \chi(X^g/{Z_G(g)})=\chi(IX/G),\ee
with $IX=\bigsqcup_{g \in G} X^g$ the inertia space as in the Introduction.

Let $H$ be the functor $F(-)$ of $\bQ$-valued algebraically constructible functions, which is covariant for all morphisms, and with a compatible cross-product. Following Ohmoto's notations, for a fixed group $A$, let $j_r(A)$ be the number of index $r$ subgroups of $A$, which is assumed to be finite for all $r$. In the notations of Theorem \ref{main2}, let 
$$b_r:=j_r(A)\cdot 1_Z \in F(Z).$$
Then  
$${\underline{b}}^{(\sigma)}={\rm Ind}^{\Sigma_n}_{Z_{\Sigma_n}(\sigma)}\left( \boxtimes_r (j_r(A)\cdot 1_Z)^{\boxtimes k_r} \right) \in  F^{\Sigma_n}(Z^n),$$
and the element $$\mathds{1}^{A}_{Z^n;\Sigma_n}:=\sum_{(\sigma)\in (\Sigma_n)_*} \underline{b}^{(\sigma)}
\in F^{\Sigma_n}(Z^n)$$ appearing on the left hand side of Equation (\ref{m2})
is a delocalized version of Ohmoto's canonical constructible function 
$\mathds{1}^{(A)}_{Z^n;\Sigma_n}$ of \cite{Oh2}[Defn.2.2],
in the sense that $$av_n \left( \mathds{1}^{A}_{Z^n;\Sigma_n}\right)=\mathds{1}^{(A)}_{Z^n;\Sigma_n}  \in F^{\Sigma_n}_{id}(Z^n).$$ This identification follows from \cite{Oh2}[Lem.3.4]. (Note that Ohmoto's product in loc.cit. corresponds to our twisted product $\widetilde{\odot}$.) 

\medskip

Let us illustrate what these distinguished constructible functions are in the case $A=\bZ$ and, resp., $\bZ^2$ (for other examples, see \cite{Oh2}). 
\begin{enumerate}
\item[(a)] if $A=\bZ$, then $$\mathds{1}^{\bZ}_{Z^n;\Sigma_n}=\bigoplus_{\sigma \in \Sigma_n} 1_{(Z^n)^{\sigma}} \in \big( \bigoplus_{\sigma \in \Sigma_n} F((Z^n)^{\sigma})  \big)^{\Sigma_n}$$
and $$\mathds{1}^{(\bZ)}_{Z^n;\Sigma_n}=\frac{1}{n!}\sum_{\sigma \in \Sigma_n} 1_{(Z^n)^{\sigma}} \in F(Z^n)^{\Sigma_n}$$
\item[(b)] if $A=\bZ^2$, then $$\mathds{1}^{(\bZ^2)}_{Z^n;\Sigma_n}=\frac{1}{n!}\sum_{\sigma \sigma' =\sigma'\sigma} 1_{(Z^n)^{\sigma} \cap (Z^n)^{\sigma'}} \in F(Z^n)^{\Sigma_n},$$ as shown in \cite{Oh2}. On the other hand, the combinatorics used in \cite{Oh2} only applies to $F(Z^n)^{\Sigma_n}$, but not to the delocalized theory $\big( \bigoplus_{\sigma \in \Sigma_n} F((Z^n)^{\sigma})  \big)^{\Sigma_n}$. Nevertheless, the function $\mathds{1}^{\bZ^2}_{Z^n;\Sigma_n}$ is a canonical lift of $\mathds{1}^{(\bZ^2)}_{Z^n;\Sigma_n}$ with respect to the averaging $av_n$.
\end{enumerate}

\medskip

Then our main Theorem \ref{main2} yields the following identity:
\be\label{o1} \sum_{n \geq 0} \ \mathds{1}^{A}_{Z^n;\Sigma_n}  \cdot t^n=\exp \left( \sum_{r\geq 1} \frac{ j_r(A)}{r} t^r \cdot \mathfrak{a}_r (1_Z) \right) \in \mathds{F}^{\Sigma}(Z).
\ee
By applying to (\ref{o1}) the ring homomorphism $av:(\mathds{F}^{\Sigma}(Z),\odot) \to (\mathds{F}_{id}^{\Sigma}(Z),\widetilde{\odot})$, we recover by (\ref{cr}) Ohmoto's generating series formula \cite{Oh2}[Prop.3.9]:
\be\label{o2}  \sum_{n \geq 0} \ \mathds{1}^{(A)}_{Z^n;\Sigma_n}  \cdot t^n=\exp \left( \sum_{r\geq 1} \frac{ j_r(A)}{r} t^r \cdot \Delta_{r*} (1_Z) \right).
\ee

Recall now that MacPherson's Chern class transformation (with rational coefficients)  $c_*:F(-) \to H_*(-):=H_{ev}^{BM}(-)\otimes \bQ$ commutes with proper pushforward and cross-products.
Applying $c_*$ to equation (\ref{o1}) we obtain a new generating series:
\be\label{o3} \sum_{n \geq 0} \ c_*(\mathds{1}^{A}_{Z^n;\Sigma_n})  \cdot t^n=\exp \left( \sum_{r\geq 1} \frac{ j_r(A)}{r} t^r \cdot \mathfrak{a}_r (c_*(1_Z)) \right) \in \mathds{H}^{\Sigma}(Z), 
\ee
with $c_*(\mathds{1}^{A}_{Z^n;\Sigma_n}) \in H^{\Sigma_n}(Z^n)$ a delocalized version of Ohmoto's orbifold Chern class $c_*(\mathds{1}^{(A)}_{Z^n;\Sigma_n}) \in H_{id}^{\Sigma_n}(Z^n)$. Ohmoto's formula \cite{Oh2}[Thm.1.1] for orbifold Chern classes of symmetric products is obtained by applying $c_*$ to (\ref{o2}), or equivalently, by applying $av$ to (\ref{o3}).

If $Z$ is projective, by taking the degrees $deg \ c_*(\mathds{1}^{A}_{Z^n;\Sigma_n})=deg \ c_*(\mathds{1}^{(A)}_{Z^n;\Sigma_n})$ of the above characteristic class formulae one recovers generating series for orbifold-type Euler characteristics, see \cite{Oh2} for details and examples.

Finally, note that for $A=\bZ$ (hence $j_r(A)=1$ for all $r$), we recover a special case (for $cl_*=c_*$ and $\cF=\bQ_Z$) of Theorem \ref{thint} from the Introduction, via the identification: $$c_*^{\Sigma_n}(\bQ_Z^{\boxtimes n})=\bigoplus_{\sigma \in \Sigma_n} c_*(1_{(Z^n)^{\sigma}})=c_*(\mathds{1}^{\bZ}_{Z^n;\Sigma_n}) \in \big( \bigoplus_{\sigma \in \Sigma_n} H_*((Z^n)^{\sigma})  \big)^{\Sigma_n}.$$ The corresponding degree formula for $Z$ compact is the classical Euler characteristic formula:
$$\chi(Z^n/\Sigma_n)=\frac{1}{n!}\sum_{\sigma\in \Sigma_n} \chi((Z^n)^{\sigma}).$$
Similarly, for $A=\bZ^2$ and $Z$ compact, one recovers at the degree level the orbifold Euler characteristic $\chi(Z^n,\Sigma_n)$ as in (\ref{oe}).


\br Instead of fixing the group $A$ and the coefficients  $b_r=j_r(A)\cdot 1_Z \in F(Z)$, one could also start with $b_r=1_Z$ for all $r$.  Applying the Frobenius-type homomorphism $av_F$ to the corresponding identity derived from our Theorem \ref{main2}, we recover the above results in a uniform way by specializing $p_r$, for a given group $A$, to $p_r=j_r(A)$ for all $r$.
\er

\br The delocalized equivariant homology $H^G_*(X)=H_*(IX/G)$ of the $G$-space $X$ and related invariants (e.g., Euler characteristic, Hodge numbers, etc.) can be used in two different ways:
\begin{enumerate}
\item[(a)] for the study of (equivariant) invariants of the quotient space $X/G$ (resp., the $G$-space $X$), via fixed-point contributions; this is the context studied in this paper via characteristic classes of Lefschetz-type.
\item[(b)] for studying (equivariant) orbifold-type invariants of $X/G$, defined as the corresponding (equivariant) invariants of the quotient space $IX/G$ (resp., the inertia space $IX$). A de Rham type description of the orbifold cohomology is already implicit in \cite{Ka}, but see also \cite{FPS} for a more recent version adapted to compact group actions.
Generating series for the corresponding orbifold Euler characteristic and orbifold Hirzebruch genus of symmetric products have been obtained, e.g., in \cite{QW,W,WZ,Zh}. Orbifold Chern class versions have been systematically studied in \cite{Oh2}, whereas orbifold elliptic classes are considered in relation to McKay correspondence in \cite{BL}. Note that the orbifold elliptic classes specialize to orbifold Hirzebruch as well as orbifold Todd classes, which are only implicitly studied in \cite{BL}.
\end{enumerate}

Let us finally note that in order to study these orbifold-type invariants via techniques of the present paper one has to apply the Lefschetz-type Riemann-Roch transformations to the inertia space $IX$ (as opposed to the $G$-space $X$, as in this paper). This in turn changes the underlying combinatorics, and it is only for constructible functions that the corresponding invariants can be directly deduced from our abstract generating series formula (\ref{m2}), as indicated above.

\er


\section{Generating series for (equivariant) characteristic classes}

In this section, we specialize our abstract generating series formula (\ref{m2}) in the framework of characteristic classes of singular varieties.

\subsection{Characteristic classes of Lefschetz type}\label{s3b}

For a complex quasi-projective variety $X$, with an algebraic action $G\times X \to X$ of a finite group $G$, let $\pi:X \to X':=X/G$ be the quotient map. We denote generically by $cat^G(X)$ a category of $G$-equivariant objects on $X$ in the underlying category $cat(X)$, see \cite{MS,CMSS}.  From now on, $H(X):=H_*(X)$ will be $H^{BM}_{ev}(X) \otimes R$, the even degree Borel-Moore homology of $X$ with $R$-coefficients, for $R$ a commutative $\bC$-algebra, resp. $\bQ$-algebra if $G$ is a symmetric group. Note that $H(-)$ is functorial for all proper maps, with a compatible cross-product (as used in the previous section).

\bd\label{Lt} An {\it equivariant characteristic class transformation of Lefschetz type} is a transformation
$$cl_*(-;g):K_0(cat^G(X)) \to H_{*}(X^g)$$ so that the following properties are satisfied:
\bn
\item $cl_*(-;g)$ is covariant functorial for $G$-equivariant proper maps.
\item $cl_*(-;g)$ is multiplicative under cross-product $\boxtimes$.
\item if $X$ is a point space and $cat(pt)$ is an abelian $\bC$-linear category, then the category ${\it Vect}_{\bC}(G)$ of (finite-dimensional) complex $G$-representations is a  subcategory of $cat^G(pt)$ and $cl_*(-;g)$ is a certain $g$-trace (as shall be explained later on), with $cl_*(-;g)= trace_g$ on $Rep_{\bC}(G)$. 
\item  if $G$ acts trivially on $X$ and $cat(X)$ is an abelian $\bC$-linear category, then $$K_0(cat^G(X)) \simeq K_0(cat(X)) \otimes Rep_{\bC}(G)$$ via the Schur functor decomposition as in (\ref{new4}), and 
\be\label{trivact} cl_*(-;g)=cl_*(-) \otimes trace_g,
\ee
with $cl_*(-)$ the corresponding non-equivariant characteristic class, as explained below. If $G=\Sigma_n$ is a symmetric group, it is enough to assume that $cat(X)$ is an abelian $\bQ$-linear category, with the category ${\it Vect}_{\bQ}(\Sigma_n)$ of rational $\Sigma_n$-representations a subcategory of $cat^{\Sigma_n}(pt)$.
\en
\ed

\br  For a subgroup $K$ of $G$, with $g \in K$, we assume that such a transformation $cl_*(-;g)$ of Lefschetz type commutes with the  restriction functor ${\rm Res}^G_K$.  Then $cl_*(-;g)$ depends only on the action of the cyclic subgroup generated by $g$. In particular, if $g=id_G$ is the identity of $G$, we can take $K$ the identity subgroup $\{id_G\}$ with  ${\rm Res}^G_{\{id_G\}}$ the forgetful functor $${\fo}: K_0(cat^G(X)) \to K_0(cat(X)),$$ so that $cl_*(-;id_G)=cl_*(-)$ fits with a corresponding non-equivariant characteristic class.
\er

\br\label{pairing} 
The above assumptions about cross-product and restriction functors can be used to define a pairing:
$${\it Vect}_{\bC}(G) \times cat^G(X) \overset{\otimes}{\lra} cat^G(X)$$
by $$cat^G(pt) \times cat^G(X) \overset{\boxtimes}{\lra} cat^{G\times G}(pt \times X) \overset{{\rm Res}}{\lra} cat^G(X),$$
with $pt \times X \cong X$, and ${\rm Res}$ denoting the restriction functor for the diagonal subgroup $G \hookrightarrow G \times G$. This induces a pairing 
$$Rep_{\bC}(G) \times K_0(cat^G(X)) \overset{\otimes}{\lra} K_0(cat^G(X))$$
on the corresponding Grothendieck groups, such that 
\be\label{npairing} cl_*(V \otimes \cF;g)=trace_{g}(V) \cdot cl_*(\cF;g)\ee
for $V$ a $G$-representation and $\cF \in cat^G(X)$. If $G$ is the symmetric group, then the above holds also for rational representations.
\er

Let us give some examples of equivariant characteristic class transformations of Lefschetz type.
\bex\label{td} {\it Todd classes} \\
Let $X$ be a quasi-projective $G$-variety, and denote by $K_0(\co^G(X))$ the Grothendieck group of the abelian category $\co^G(X)$ of $G$-equivariant coherent algebraic sheaves on $X$. For each $g \in G$, the Lefschetz-Riemann-Roch transformation \cite{BFQ,Mo}
\begin{equation}
td_*(-;g):K_0(\co^G(X)) \lra H_{*}(X^g)
\end{equation}
is of Lefschetz type with $R:=\bC$ (resp. $R:=\bQ$ if $G$ is a symmetric group). Moreover, $td_*(-;id_G)$ is the complexified (non-equivariant) Todd class transformation $td_*$ of Baum-Fulton-MacPherson \cite{BFM}.
Over a point space, the transformation $td_*(-;g)$ reduces to the $g$-trace on the corresponding $G$-equivariant vector space. In particular, if $X$ is projective, by pushing down to a point we recover the equivariant holomorphic Euler characteristic, i.e., for $\cF\in \co^G(X)$ the following degree formula holds:
$$\chi_a(X,\cF;g):=\sum_i (-1)^i trace\left(g \vert H^i(X;\cF) \right) = \int_{[X^g]} td_*([\cF];g).$$
\eex

\bex\label{ch} {\it Chern classes}\\
Let $K_0({\cons}^G(X))$ be the Grothendieck group of the abelian category ${\cons}^G(X)$ of algebraically constructible $G$-equivariant sheaves of complex vector spaces on $X$. Then the localized Chern class transformation \cite{Sch1}[Ex.1.3.2]
$$c_*(-;g):=c_*(tr_g(-\vert_{X^g})):K_0({\cons}^G(X)) \lra H_{*}(X^g)$$
is of Lefschetz type  with $R:=\bC$ (resp. $R:=\bQ$ if $G$ is a symmetric group). Here $c_*(-)$ is the Chern-MacPherson class transformation \cite{MP}, and $tr_g(-\vert_{X^g})$ is the group homomorphism which, for $\cF \in {\cons}^G(X)$, assigns to $[\cF] \in K_0({\cons}^G(X))$ the constructible function on $X^g$ defined by $$x \mapsto trace(g\vert \cF_x).$$ (Note that for $x \in X^g$, $g$ acts on the finite dimensional stalk $\cF_x$ for a constructible $G$-equivariant sheaf $\cF$.) For the identity element, the transformation $c_*(-;g)$ reduces to the complexification of MacPherson's Chern class transformation. It also follows by definition that if $X$ is a point space, then $c_*(-;g)$ reduces to the $g$-trace on the corresponding $G$-equivariant vector space.  In particular, if $X$ is projective, by pushing down to a point we recover the equivariant Euler characteristic, i.e., for $\cF\in {\cons}^G(X)$  the following degree formula holds:
$$\chi(X,\cF;g):=\sum_i (-1)^i trace\left(g \vert H^i(X;\cF) \right) = \int_{[X^g]} c_*([\cF];g).$$
\eex

\bex {\it (un-normalized) Atiyah-Singer  classes -- mixed Hodge module version} \\
Let $K_0 (\MHM^G(X))$ be the Grothendieck group of $G$-equivariant (algebraic) mixed Hodge modules. The (un-normalized) Atiyah-Singer class transformation of \cite{CMSS},
$$T_{y*}(-;g): K_0 (\MHM^G(X)) \lra H_{*}(X^g)$$
is of Lefschetz type with $R=\bC[y^{\pm 1}]$ (resp. $R:=\bQ[y^{\pm 1}]$ if $G$ is a symmetric group). 
In this case, $\MHM(X)$ is only a $\bQ$-linear abelian category, so for the isomorphism of Grothendieck groups in property  $(3)$ of Definition \ref{Lt} one should assume that $G$ is a symmetric group (or work with the Grothendieck group of the underlying $\bC$-linear exact category of filtered holonomic $\cD$-modules).
For the identity element of $G$, this reduces to the mixed Hodge module version of the (un-normalized) Hirzebruch class transformation of Brasselet-Sch\"urmann-Yokura \cite{BSY,Sch2}. Over a point space, $T_{y*}(-;g)$ coincides with the equivariant $\chi_y(-;g)$-genus ring homomorphism 
$$\chi_y(-;g):K_0^G(\mhs) \to \bC[y^{\pm 1}]$$
defined on the Grothendieck group of the category $G-\mhs$ of $G$-equivariant (graded) 
polarizable mixed Hodge structures, by:
$$\chi_y([H];g):=\sum_p trace\left(g\vert Gr^p_F(H \otimes \bC)\right) \cdot (-y)^p,$$
for $F^{\centerdot}$ the Hodge filtration on $H \in G-\mhs$. Here we use the identification $\MHM^G(pt) \simeq G-\mhs$ of $G$-equivariant mixed Hodge modules over a point space with $G$-equivariant (graded) polarizable mixed Hodge structures, so that the category of (finite-dimensional) rational $G$-representations is a subcategory of $G-\mhs$ (viewed as mixed Hodge structure of pure type $(0,0)$). In particular, if $X$ is projective and $\cM \in \MHM^G(X)$, by pushing down to a point we recover the equivariant twisted Hodge genus, i.e., the following degree formula holds (see \cite{CMSS}[Prop.4.7] for $\cM$  the intersection cohomology mixed Hodge module):
$$\chi_y(X, \cM;g):=\sum_{i,p} (-1)^i trace\left(g \vert Gr^p_F H^i(X;\cM) \otimes \bC \right)\cdot (-y)^p =\int_{[X^g]} T_{y*}(\cM;g).$$
\eex

\br The motivic version of the Atiyah-Singer class transformation, as mentioned in the introduction, can be deduced from the corresponding mixed Hodge module version through 
 the natural transformation (e.g., see \cite{CMSS})
$$\chi^G_{\rm Hdg}: K^G_0 (var/X) \to K_0 (\MHM^G(X))$$ mapping $[id_X]$ to the class of the constant Hodge module (complex) $[\bQ^H_X]$. This transformation commutes with push downs, cross-products and restriction functors. Then the (un-normalized) motivic Atiyah-Singer class transformation is the composition
$$T_{y*}(-;g):K^G_0 (var/X) \overset{\chi^G_{\rm Hdg}}{\lra} K_0 (\MHM^G(X)) \overset{T_{y*}(-;g)}{\lra} H_{*}(X^g).$$
In particular, the Atiyah-Singer class of $X$ is defined as 
$$T_{y*}(X;g):=T_{y*}([id_X];g)=T_{y*}([\bQ^H_X];g).$$
\er

\bigskip

Let us explain some of the above examples in the simplest situation when $X$ is  {\it smooth} (see \cite{CMSS} for complete details). Then $X^g$ is also smooth, and we denote by $T_{X^g}$ and $N_{X^g}$ its tangent and resp. normal bundle in $X$. In this case, the homological Leschetz-type transformations correspond under Poincar\'e duality  (and for suitable ``smooth'' coefficients in $cat^G(X)$) to similar cohomological transformations, as explained below. 

\medskip

\noindent{\bf Todd classes.} Let $K_G^0(X)$ be the Grothendieck group of algebraic $G$-vector bundles, and note that the natural map $K_G^0(X) \to K_0(\co^G(X))$ is an isomorphism. Let $ch^*$ and $td^*$ denote the Chern character and resp. the Todd class in cohomology. The Lefschetz-Riemann-Roch transformation is then given by: for  $V$ an algebraic $G$-vector bundle on $X$,
\be\label{AS1}
td_*(V;g)=ch^*(g)(V|_{X^g}) \cap \left(\frac{td^*(T_{X^g})}{ch^*(g)(\Lambda_{-1}N^*_{X^g})} \cap [X^g] \right).
\ee
Here, $N^*_{X^g}$ denotes the dual of the normal bundle of $X^g$, and for a vector bundle $E$ we let $\Lambda_{-1}(E):=\sum_i (-1)^i \Lambda^i E$. Moreover, the equivariant Chern character $ch^*(g)(-):K_G^0(X^g) \to H^*(X^g)$ is defined as 
follows: for $W \in K_G^0(X^g)$ we let
$$ch^*(g)(W):=\sum_{\chi} \chi(g) \cdot ch^*(W_{\chi}),$$
for $W \simeq \oplus_{\chi} W_{\chi}$ the (finite) decomposition of $W$ into sub-bundles $W_{\chi}$ on which $g$ acts by a character $\chi: \langle g \rangle  \to \bC^*$. Note that $ch^*(g)(V|_{X^g})$ is just the complexified Chern character of $L(g)(V)$, with $L(g)$ the Lefschetz-type transformation of Section \ref{loc1}.
If $X$ is projective, by taking degrees in formula (\ref{AS1}) we obtain the Atiyah-Singer holomorphic Lefschetz formula from \cite{AS,HZ}.

\medskip

\noindent{\bf Atiyah-Singer classes (mixed Hodge module version).}   Let $X$ be smooth with an algebraic $G$-action, together with a $G$-equivariant ``good'' variation $\cL$ of rational mixed Hodge structures (i.e., graded polarizable, admissible and with quasi-unipotent monodromy at infinity). This corresponds to a (shifted) smooth $G$-equivariant mixed Hodge module. Let $\mathcal{V}:=\cL \otimes_{\bQ} \mathcal{O}_X$ be the flat $G$-equivariant bundle whose sheaf of
horizontal sections is $\cL \otimes \bC$. The bundle $\cV$ comes
equipped with a decreasing (Hodge) filtration (compatible with the $G$-action) by holomorphic
sub-bundles $\cF^p\cV$. Note that since we work with a ``good'' variation, each $\cF^p\cV$ underlies (by GAGA) a unique complex algebraic $G$-vector bundle. Let $$\chi_y(\cV):=\sum_p \left[{\rm Gr}^{p}_{\mathcal{F}} \cV \right]
\cdot (-y)^{p} \in K^0_G(X)[y,y^{-1}]$$
be the $\chi_y$-characteristic of $\cV$. Then:
\be{{T_y}}_*(X,\cL;g)
=ch^*(g)(\chi_y(\cV) \vert_{X^g}) \cap T_{y*}(X;g),\ee
with $T_{y*}(X;g):=\sum_{i \geq 0}td_*([\Omega^i_X];g) \cdot y^i$ the Atiyah-Singer class of $X$.

\medskip

Our generating series results for these characteristic class transformations, as discussed in the next sections will, however, be valid for any quasi-projective complex variety $X$ (possibly singular) and all coefficients in $cat^G(X)$.



\subsection{Delocalized equivariant characteristic classes} \label{dlcl}
Let $X$ be a (possibly singular) quasi-projective variety acted upon by a finite group $G$ of algebraic automorphisms.

From now on, we use the symbol $cl_*$ to denote any of the characteristic classes $c_*$, $td_*$ and $T_{y*}$, respectively, with their corresponding equivariant versions of Lefschetz type, $cl_*(-;g):K_0(cat^G(X)) \to H_*(X^g)$, as discussed in the previous section. 
\bd For any of the Lefschetz-type characteristic class transformations $cl_*(-;g)$ considered above, we define a  corresponding  $G$-equivariant  class transformation
$$cl^G_{*}: K_0(cat^G(X)) \to  H^{G}_{*}(X)$$
 by: $$cl^G_{*}(-):=\bigoplus_{g \in G} cl_*(-;g)=\bigoplus_{(g)} {\rm Ind}_{Z_G(g)}^G \left({cl}_*(-;g)\right) \in \left( \bigoplus_{g \in G} \ H_*(X^g) \right)^G,$$
 with induction as in Remark \ref{expl}.
\ed
Note that the $G$-invariance of the class $cl_*^G(-)$ is a consequence of {\it conjugacy invariance} of  ${cl}_*(-;g)$, proved in \cite{CMSS}[Sect.5.3.]. This also explains the equality of the two descriptions of $cl^G_{*}(-)$.

\medskip

The above transformation $cl^G_*(-)$ has the same properties as the Lefschetz-type transformations $cl_*(-;g)$, e.g., functoriality for proper push-downs, restrictions to subgroups, and multiplicativity for exterior products. 

\medskip

If $X$ is projective, then by pushing $cl^G_*(-)$ down to a point, the degree
$$
deg \left(cl^G_*(-)\right)=tr_G(X,-):=tr_G (H^*(X,-)) \in C(G) \otimes R
$$
is the character $tr_G$ of the corresponding virtual cohomology representation $$\sum_i (-1)^i [H^i(X,-)] \in Rep_{\bC}(G) \otimes R  \simeq C(G) \otimes R.$$
In addition, if $cl^G_*=T^G_{y*}$ is the $G$-equivariant Hirzebruch class, then $$deg \left(T^G_{y*}([\cM])\right)=tr_G(X,\cM):=\sum_p tr_G (Gr^p_FH^*(X,\cM)) (-y)^p,$$ for $\cM \in \MHM(X)$.

\medskip

If $G$ acts trivially on $X$, then there is a functor 
$$[-]^G:K_0(cat^G(X)) \to K_0(cat(X))$$ 
defined by taking $G$-invariants,
which for any of the $\bQ$-linear abelian categories $\co(X)$, $\cons(X)$ and $\MHM(X)$ is
induced from the exact projector
$$(-)^G:=\frac{1}{|G|} \sum_{g \in G} \mu_g: cat^G(X) \to cat(X).$$
Here $\mu_g:\cF \to g_*\cF$ (with $g \in G$) is the isomorphism of the $G$-action on $\cF \in cat^G(X)$, see \cite{CMSS,MS} for details. 
\br\label{Schur} For a given $G$-representation $V$, by using the pairing of Remark \ref{pairing} one can define {\it Schur functors} $S_V:cat^G(X) \to cat(X)$ by $S_V(\cF):=(V \otimes \cF)^G$. Here we assume that $G$ acts trivially on $X$.  This notion agrees with the abstract categorical definition from \cite{D,He}.
\er
For the Grothendieck group $K_0^G(var/X)$ of $G$-varieties over $X$,  the functor $$[-]^G:K_0^G(var/X) \to K_0(var/X)$$ is given by: $[Y \to X] \mapsto [Y/G \to X]$. This is a well-defined functor since our equivariant Grothendieck group $K_0^G(var/X)$ from \cite{CMSS} only uses the scissor relation.
Note that the transformation $$\chi^G_{\rm Hdg}: K^G_0 (var/X) \to K_0 (\MHM^G(X))$$ relating the two versions of the Atiyah-Singer transformation, commutes with $[-]^G$ since, for $\pi:Y \to Y/G$ denoting the quotient map, we have that (cf. \cite{CMSS}[Lem.5.3], but see also \cite{MSS}[Rem.2.4]) $$\bQ^H_{Y/G}=(\pi_*\bQ^H_Y)^{G} \in D^b\MHM(Y/G).$$

Then if $G$ acts trivially on $X$, the following {\it averaging property} holds by the definition of the projector $(-)^G$ together with (\ref{trivact}) (compare also with \cite{CMSS}[Sect.5.3], \cite{CMSSY}[Sect.3]): 
\begin{equation}\label{mainav}\begin{CD}
K_0(cat^G(X)) @> cl^G_{*} >> H^G_{*}(X) \cong H_*(X) \otimes C(G)  \\
@V [-]^G\ VV @VV {\frac{1}{\vert G \vert}\sum_{g \in G} ev_g} V \\
K_0(cat(X)) @> cl_{*} >> H_{*}(X) ,
\end{CD}\end{equation}
where $ev_g$ is the evaluation at $g \in G$ of class functions on $G$.


\subsection{Proof of the main Theorem \ref{thint} and its applications.}\label{gensersp}
In this section, we explain how to deduce our main Theorem \ref{thint} for equivariant characteristic classes of external products of varieties from the abstract generating series formula (\ref{m2}). We also explain the various specializations of Theorem \ref{thint}, as formulated in the Introduction.

\medskip

Let $Z$ be a quasi-projective variety, with a given object $\cF \in cat(Z)$ as in the Introduction. We use as before the symbol $$cl_*:K_0(cat(Z)) \to H_*(Z)=H^{BM}_{2*}(Z) \otimes R$$ to denote any of the following characteristic classes
\begin{itemize}  
\item $td_*:K_0({\co}(Z)) \to H^{BM}_{2*}(Z) \otimes \bQ$, 
\item $c_*:K_0({\cons}(Z)) \to H^{BM}_{2*}(Z) \otimes \bQ$, 
\item $T_{y*}:K_0(var/Z) \to H^{BM}_{2*}(Z) \otimes \bQ[y]$,
\item  $T_{y*}:K_0(\MHM(Z)) \to H^{BM}_{2*}(Z) \otimes \bQ[y^{\pm 1}]$, 
\end{itemize} 
with their corresponding localized and delocalized equivariant versions $cl_*(-;g)$ and resp. $cl_*^G(-)$, as discussed in the previous two subsections.

\medskip

The following properties will allow us to further specialize our main generating series formula (\ref{m2}) in the context of characteristic classes.

It follows from  \cite{CMSSY} that $cl_*(-;g)$ satisfies the following {\it multiplicativity property} (see \cite{CMSSY}[Lemma 3.2, Lemma 3.5, Lemma 3.9]): 
\bl\label{mult} (Multiplicativity)\newline
If $\sigma \in \Sigma_n$ has cycle-type $(k_1, k_2, \cdots)$, i.e., $k_r$ is the number of $r$-cycles in $\sigma$ and $n=\sum_r r \cdot k_r$, then:
\be cl_*(\cF^{\boxtimes n};\sigma)=\boxtimes_r \left(cl_*(\cF^{\boxtimes r}; \sigma_r ) \right)^{k_r} \in H_*((Z^n)^{\sigma})^{Z_{\Sigma_n}(\sigma)} \subset H_*((Z^n)^{\sigma}),
\ee
with $\sigma_r$ denoting an $r$-cycle in $\Sigma_r$.
\el

Moreover, the following {\it localization} result holds, see \cite{CMSSY}[Lemma 3.3, Lemma 3.6, Lemma 3.10]:
\bl\label{loc} (Localization)\newline 
Under the identification $(Z^r)^{\sigma_r}\simeq Z$, the following holds:
\be\label{lf}
cl_*(\cF^{\boxtimes r}; \sigma_r)=\Psi_r cl_*(\cF), 
\ee
where $\Psi_r$ denotes the homological Adams operation defined by
$$\Psi_r=
\begin{cases}
id & {\rm if} \ cl_*=c_*\\
\cdot \frac{1}{r^i} \ {\rm on} \ H^{BM}_{2i}(Z) \otimes \bQ & {\rm if} \ cl_*=td_* \\
\cdot \frac{1}{r^i} \ {\rm on} \ H^{BM}_{2i}(Z) \otimes \bQ, {\rm and} \ y \mapsto y^r & {\rm if} \ cl_*=T_{-y*}.
\end{cases}
$$
\el

\bigskip

We can now explain how to use the abstract generating series formula (\ref{m2}) to derive Theorem \ref{thint} from the Introduction.

\begin{proof} [Proof of Theorem \ref{thint}]
For any $r \geq 1$, let $$b_r:=cl_*(\cF^{\boxtimes r}; \sigma_r) \in H_*(Z).$$
By multiplicativity (cf. Lemma \ref{mult}) and conjugacy invariance of $cl_*(-;\sigma)$, it follows that  
$$cl_*^{\Sigma_n} (\cF^{\boxtimes n})= \sum_{(\sigma) \in (\Sigma_n)_*} {\rm Ind}^{\Sigma_n}_{Z_{\Sigma_n(\sigma)}}\left( \boxtimes_r cl_*(\cF^{\boxtimes r};\sigma_r)^{\boxtimes k_r} \right)= \sum_{(\sigma) \in (\Sigma_n)_*} \underline{b}^{(\sigma)}.$$
Then (\ref{m2i}) follows from our main formula (\ref{m2}) together with the localization formula (\ref{lf}).
\end{proof}

Let us now apply $\pi_*:=\oplus_n {\pi_n}_*$ to formula (\ref{m2i}). Then, by using functoriality and the corresponding $\bQ$-algebra homomorphism $\pi_*:\mathds{H}^{\Sigma}_*(Z) \to \mathds{PH}_*^{\Sigma}(Z)$ of (\ref{two}), we obtain by the first identity of (\ref{cr}) a proof of Corollary \ref{cormaini} from the Introduction.

\medskip

The averaging property  (\ref{mainav}) together with the $\bQ$-algebra evaluation homomorphism $\mathds{PH}_*(Z) \otimes \bQ[p_i , i \geq 1] \to \mathds{PH}_*(Z)$, $p_i \mapsto 1$ (for all $i$), as in (\ref{two}), now yields a proof of Corollary \ref{c1.6} from the Introduction, which recovers the main result of \cite{CMSSY}.

\medskip

We want to emphasize that Corollary \ref{cormaini} has other important applications, as already mentioned in the Introduction. For example, it specializes to Corollary \ref{c1.7} by using the $\bQ$-algebra evaluation homomorphism $$\mathds{PH}_*(Z) \otimes \bQ[p_i , i \geq 1] \to \mathds{PH}_*(Z), \ \ p_i \mapsto sign(\sigma_i)=(-1)^{i-1}$$ (for all $i$), together with the commutative diagram:
\begin{equation}\label{mainavsp}\begin{CD}
K_0(cat^{\Sigma_n}(\Zs)) @> cl^{\Sigma_n}_{*} >> H^{\Sigma_n}_{*}(\Zs) \cong H_*(\Zs) \otimes C(\Sigma_n)  \\
@V [-]^{sign-\Sigma_n}\ VV @VV {\frac{1}{n!}\sum_{\sigma \in \Sigma_n} sign(\sigma) \cdot ev_{\sigma}} V \\
K_0(cat(\Zs)) @> cl_{*} >> H_{*}(\Zs) .
\end{CD}\end{equation}

\bigskip

Finally, formula (\ref{new2}) from the Introduction follows by combining the multiplicativity of (\ref{npairing}), together with the identification of the coefficient of $t^n$ in the explicit expansion (as in the proof of Theorem \ref{main2}) of the exponential on the right-hand side of (\ref{ci}), that is,
\be cl_*^{\Sigma_n} (\pi_{n*}\cF^{\boxtimes n})=\sum_{{\lambda=(k_1,k_2, \cdots) \dashv \ n} } \frac{p_{\lambda}}{z_{\lambda}}  \cdot \bigodot_{r \geq 1} \left( d_{r*} (\psi_r(cl_*(\cF)) \right)^{k_r}.\ee

\br\label{5.9}
All these results and arguments also apply to a bounded complex $\cF$ in $D^b_{coh}(Z)$, resp. $D^b_c(Z)$, with coherent, respectively constructible cohomology, as well as to bounded complexes  (such as $\bQ^H_Z$) in  $D^b(\MHM(Z))$. Then $\cF^{\boxtimes n}$ and $\pi_{n*}\cF^{\boxtimes n}$ become {\it weakly equivariant $\Sigma_n$-complexes} (as in \cite{CMSS}[Appendix]), which still have well-defined Grothendieck classes $[\cF^{\boxtimes n}] \in K_0(cat^{\Sigma_n}(Z^n))$, resp.,  $[\pi_{n*}\cF^{\boxtimes n}]\in K_0(cat^{\Sigma_n}(\Zs))$. Moreover, the definition of the symmetric and resp. alternating power objects via the projector $[-]^{\Sigma_n}$ and resp. $[-]^{sign-\Sigma_n}$  still works as above since the corresponding derived categories are $\bQ$-linear Karoubian (see \cite{CMSSY, MS} for more details).  
\er

\

\end{document}